\documentclass[sigconf]{acmart}
\usepackage[utf8]{inputenc}
\usepackage{graphics} 
\usepackage{xcolor}
\usepackage{graphicx}
\usepackage{caption}
\usepackage{amsmath,amsfonts}
\usepackage{bbm, bm}
\usepackage{subcaption}
\usepackage{listings}
\usepackage{hyperref}
\usepackage{pifont}
\usepackage{tabularx}
\usepackage{algorithmicx}
\usepackage[noend]{algpseudocode}
\usepackage{algorithm}
\usepackage{tablefootnote}
\usepackage{float}
\usepackage{nicematrix}
\usepackage{tikz}
\usepackage{booktabs}
\usetikzlibrary{arrows}
\usetikzlibrary{positioning,matrix,calc}

\newcommand{\coo}{\ensuremath{\mathrm{CO_2}}} 
\newcommand{\trianglebullet}{$\mbox{\ensuremath{\rhd}}$}
\newcommand{\numsim}{$500$~} 

\newcommand{\multicell}[1]{\begin{tabular}[c]{@{}c@{}} #1 \end{tabular}}
\newcommand{\multicellleft}[1]{\begin{tabular}[l]{@{}l@{}} #1 \end{tabular}}

\ifodd 1

    \newcommand{\com}[2]{\textbf{\color{blue} (COMMENT from [#1]: #2)}}
\else

    \newcommand{\com}[2]{}
\fi

\AtBeginDocument{%
  \providecommand\BibTeX{{%
    \normalfont B\kern-0.5em{\scshape i\kern-0.25em b}\kern-0.8em\TeX}}}

\copyrightyear{2023} 
\acmYear{2023} 
\setcopyright{acmlicensed}\acmConference[e-Energy '23]{The 14th ACM
International Conference on Future Energy Systems}{June 20--23,
2023}{Orlando, FL, USA}
\acmBooktitle{The 14th ACM International Conference on Future Energy
Systems (e-Energy '23), June 20--23, 2023, Orlando, FL, USA}
\acmPrice{15.00}
\acmDOI{10.1145/3575813.3595193}
\acmISBN{979-8-4007-0032-3/23/06}

\begin{document}

\settopmatter{printfolios=true}

\title{Follow the Sun and Go with the Wind:\\ Carbon Footprint Optimized Timely E-Truck Transportation}

\author{Junyan Su}
\affiliation{%
  \institution{School of Data Science\\City University of Hong Kong}
  \city{}
  \country{}
}
\author{Qiulin Lin}
\affiliation{%
  \institution{School of Data Science\\City University of Hong Kong}
  \city{}
  \country{}
}

\author{Minghua Chen}
\thanks{Corresponding author: Minghua Chen.}
\affiliation{%
  \institution{School of Data Science\\City University of Hong Kong}
  \city{}
  \country{}
}


\begin{abstract}
We study the carbon footprint optimization (CFO) of a heavy-duty e-truck traveling from an origin to a destination across a national highway network subject to a hard deadline, by optimizing path planning, speed planning, and intermediary charging planning.  Such a CFO problem is essential for carbon-friendly e-truck operations. However, it is notoriously challenging to solve due to (i) the hard deadline constraint, (ii) positive battery state-of-charge constraints, (iii) non-convex carbon footprint objective, and (iv) enormous geographical and temporal charging options with diverse carbon intensity. Indeed, we show that the CFO problem is NP-hard. As a key contribution, we show that under practical settings it is equivalent to finding a generalized restricted shortest path on a stage-expanded graph, which extends the original transportation graph to model charging options. Compared to alternative approaches, our formulation incurs low model complexity and reveals a problem structure useful for algorithm design. We exploit the insights to develop an efficient dual-subgradient algorithm that always converges. As another major contribution, we prove that (i) each iteration only incurs polynomial-time complexity, albeit it requires solving an integer charging planning problem optimally, and (ii) the algorithm generates optimal results if a condition is met and solutions with bounded optimality loss otherwise. Extensive simulations based on real-world traces show that our scheme reduces up to {28\%}  carbon footprint compared to baseline alternatives. The results also demonstrate that e-truck reduces 56\% carbon footprint than internal combustion engine trucks.
\end{abstract}
\maketitle

\begin{figure}[!tb]
	\centering
	\includegraphics[width=\linewidth]{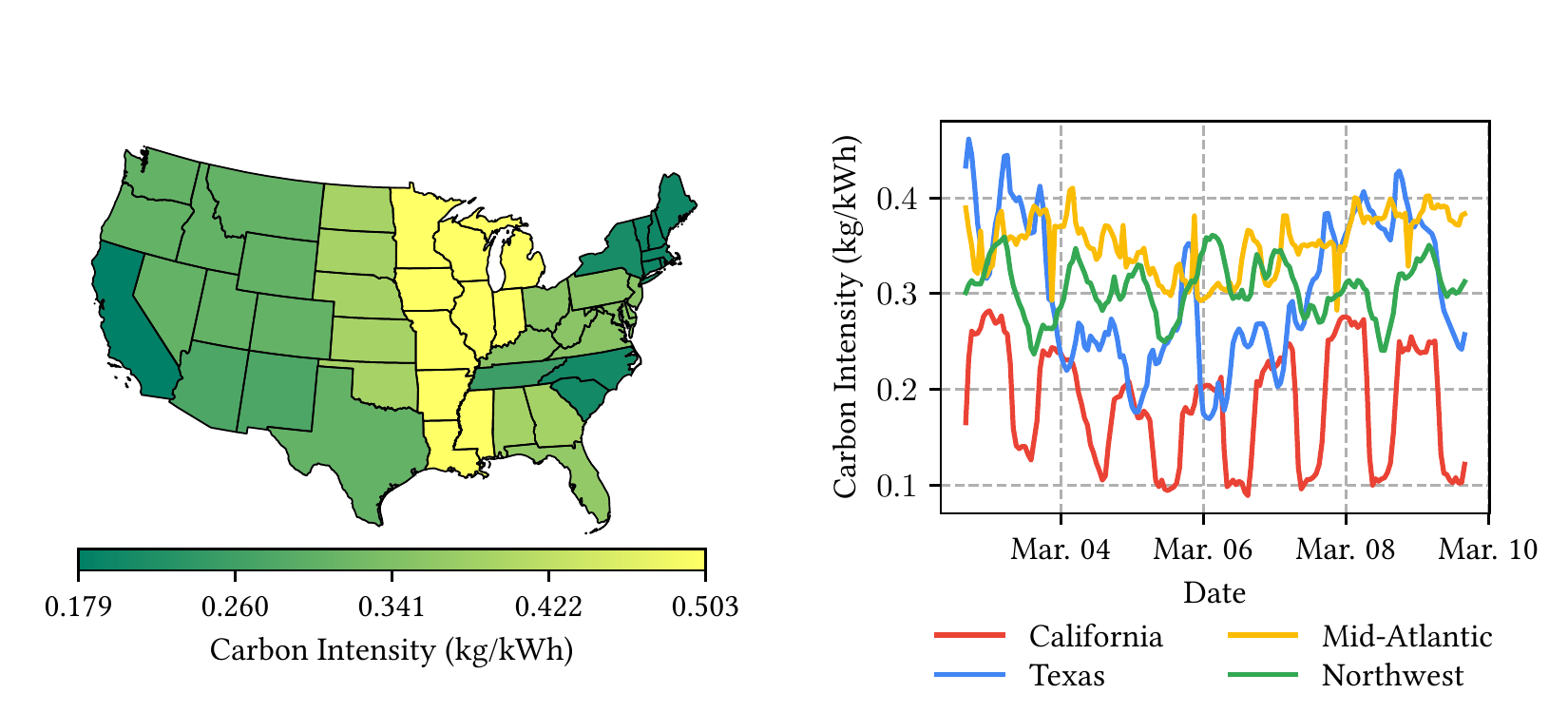} 
    \vspace{-4mm}
	\caption{Spatial-temporal diversity of the carbon intensity~(kg/kWh) in the US~\cite{energyinformationadministrationeiaRealtime}. The whole country is divided to 13 regions. Each region includes one or multiple states. Carbon intensity is only provided for individual regions~\cite{energyinformationadministrationeiaRealtime}, thus one uniform color for all the states in a region.
    }
	\label{fig:price-distri}
\end{figure}

\begin{table}[!tb]
    \centering
    \resizebox{\linewidth}{!}{
    \begin{tabular}{lllll}
        \toprule
        & Coal & Natural gas & Petroleum & Renewable \\
        \midrule
        Carbon Intensity (kg/kWh) & 1.02 & 0.39 & 0.91 & 0 \\
        Total emission (tons) & 7.86 $\times 10^8$ & 6.35 $\times 10^8$& 1.6 $\times 10^7$ & 0\\
        Electricity (kWh) & 7.73 $\times 10^{11}$ & 1.62$\times 10^{12}$ & 1.75 $\times 10^{10}$ & $7.92\times 10^{11}$  \\
        \bottomrule
    \end{tabular}
    } 
    \vspace{1mm}
    \caption{U.S electricity generation and resulting \coo\ emissions by major type of source in 2020~\cite[Tab.~7.2a, Tab.~11.6]{eia2021monthly}}\label{tab:fuel-co2}
    \vspace{-6mm}
\end{table}

\section{Introduction}

In the United States, heavy-duty trucks deliver more than 10 billion tons of freight in 2021, representing 72.2\% of total domestic tonnage shipped. 80.4\% of the nation’s freight bills are from trucking which amounts to 875 billion US dollars~\cite{truckData2021}.
Despite their importance for the economy, heavy-duty trucks are a major source of Carbon Dioxide (\coo) emission, which is tied to climate change. With only $4\%$ of the total vehicle population, heavy-duty trucks produce  $25\%$ of the \coo\ emitted in the transportation sector, accounting for $8.8\%$ of the total carbon emissions in the US~\cite{davisTransportation2021}. To fight climate change, it is critical to de-carbonize the trucking industry.

\begin{table*}[!tb]
    \centering
    
\newcommand{\cmark}{\ding{51}}
\newcommand{\xmark}{\ding{55}}%

\begin{tabular}{lccccll}
    \toprule
    & \multicell{Speed Planning} & \multicell{Path Planning} & \multicell{Deadline Constraint} & \multicell{Charging Planning} & \multicell{Truck Type} & Objective  \\ \midrule
    RSP~\cite{hassin1992approximation,lorenz2001simple,juttner2001lagrange} & \xmark & \cmark & \cmark & N/A & ICE & Path cost \\
    E2T2~\cite{dengEnergyefficient2017,liuEnergyEfficient2018,xuRide2019,zhou2020minimizing} & \cmark & \cmark & \cmark & N/A & ICE & Energy/Emission \\
    many, e.g.,~\cite{hellstrom2009look} & \cmark & \xmark & \xmark & N/A & ICE & Energy \\
    many, e.g.,~\cite{strehler2017energy} & \cmark & \cmark & \xmark & \cmark & Electric/Hybrid & Time/Energy \\
    many, e.g.,~\cite{zhang2022optimal} & \cmark & \xmark & \cmark & \xmark & Electric & Energy \\
    \textbf{This work} & \cmark & \cmark & \cmark & \cmark  & Electric & Carbon footprint \\ \bottomrule
\end{tabular}

    \vspace{1mm}
    \caption{Model comparison between this work and the existing work for truck operations. 
    }\label{tab:related-truck}
    \vspace{-0mm}
\end{table*}

Electrification of heavy-duty trucks is a promising endeavor with multi-facet benefits, including lower noise disturbance, better acceleration performance, better energy efficiency, and most notably, zero tailpipe emissions~\cite{bibra2022global} with great potential to de-carbonize the trucking industry.
However, while e-trucks do not generate direct \coo~emission during the operation, charging e-trucks incurs positive \emph{carbon footprint}, i.e., the \coo~emission during the production of electricity to be consumed by the e-truck. It is thus essential to evaluate the carbon footprint of e-trucks in contrast to the tailpipe emission for ICE trucks~\cite{zhou2020minimizing}.
{
}  

In this paper, we study the carbon footprint optimization (CFO) of a heavy-duty e-truck traveling from an origin to a destination across a national highway network subject to a hard deadline, by exploring the complete design space of path planning, speed planning, and intermediary charging planning.  Such a CFO problem is central to maximizing the environmental benefit brought by e-trucks. However, it is notoriously challenging to solve due to the following.

\emph{The hard deadline constraint.} Transportation deadline is a common requirement in the trucking industry. In fact, timely delivery is not only necessary for perishable goods such as food~\cite{ashby1987protecting}, it is also commonly adapted in service-level agreements to guarantee delivery delay ~\cite{amazon}. Moreover, transportation tasks on online platforms (e.g., Uber Freight~\cite{uber}) are often associated with pickup/delivery time requirements, from which an end-to-end delivery deadline can be computed. Meanwhile, very often, considering deadline constraints can turn an easy problem, e.g., finding the fastest path between the origin and a destination, into an NP-hard counterpart~\cite{dengEnergyefficient2017}.

\emph{Positive battery State of Charge (SoC) constraints.}
Range anxiety is a common issue for e-trucks due to their limited capacity. Careful planning to prevent battery depletion for an e-truck is thus necessary during the transportation task. However, this requirement adds strong couplings to the traveling decisions between successive road segments, thus largely complicating the problem. 


\emph{Non-convex carbon footprint objective.}
We consider an e-truck incurs the carbon footprint objective during intermediary charging.
End-use electricity often has a mix of different sources (e.g., coal, natural gas, petroleum, and renewable). Different sources have different carbon intensities (i.e., carbon footprint per unit electricity, unit: kg/kWh); see Tab.~\ref{tab:fuel-co2} for examples.
Among those sources, renewable have notable zero carbon intensity and thus is an eco-friendly type of source.
Meanwhile, renewable generation has a large variation in time and locations, resulting in spatial-temporal diversity of carbon intensity of electricity; see Fig.~\ref{fig:price-distri} for an illustration.
The diverse carbon intensity, along with the nonlinear charging characteristics of the battery, makes the carbon footprint objective non-convex and challenging to deal with.


\emph{Enormous geographical and temporal charging options.}
Due to the fluctuating nature of the renewable, carbon intensity varies largely in both location and time. 
Jointly with path planning and speed planning, charging planning determines where and when to charge, how much to charge with diverse carbon intensity in different locations and time.
The design space is combinatorial in nature, with enormous options to consider, making the CFO problem uniquely challenging to solve.

Indeed, the CFO problem is more complicated and challenging than those in related studies for ICE truck and EV driving optimization and existing approaches are not directly applicable; see Sec.~\ref{sec:related} for more discussions. 
In this paper, we carry out a comprehensive study of the problem and make the following \textbf{contributions}: 

\trianglebullet 
We conduct the first study on \emph{carbon footprint optimization} (CFO).
We show that the problem is NP-hard. 
We then show that under practical settings, the CFO problem is equivalent to finding a generalized restricted shortest path on a stage-expanded graph, which extends the original transportation graph to incorporate charging options. 
Compared to alternative approaches, our novel formulation incurs low model complexity and reveals a problem structure useful for algorithm design.


\trianglebullet 
We then exploit the structured insights from the formulation to develop an efficient dual-subgradient algorithm that always converges. 
We show that in our approach, each iteration only incurs polynomial-time complexity in graph size, albeit it requires solving an integer charging planning problem optimally.
We then further derive a condition under which our algorithm produces an optimal solution and a posterior performance bound when the condition is not satisfied.

\trianglebullet  
We carry out extensive numerical experiments using real-world traces over the US highway network to show the effectiveness of our approach in practice. 
The results show that our scheme reduces up to {28\%} carbon footprint as compared to a fastest-path alternative. 
We conduct a comparison between e-trucks and ICE trucks on their environmental impact. The results show that an e-truck achieves 56\% carbon reduction compared to its ICE counterpart.

\section{Related Work}\label{sec:related}

\textbf{Energy-Efficient Timely Transportation (E2T2).} E2T2 aims to find a path and its associated speed profile from an origin to a destination with minimum energy consumption, while satisfying a hard deadline~\cite{dengEnergyefficient2017,liuEnergyEfficient2018,xuRide2019}. It generalizes the NP-hard Restricted Shortest Path (RSP)~\cite{hassin1992approximation,lorenz2001simple,juttner2001lagrange} problem by introducing an additional design space of speed planing. 
E2T2 problem is first proposed in~\cite{dengEnergyefficient2017} and also shown to be NP-hard. E2T2 is then generalized to the scenario with multiple pickup/delivery locations with time windows in~\cite{liuEnergyEfficient2018}. Another extension of E2T2 has been studied in~\cite{xuRide2019}, where the driver can strategically wait at the rest area for better traffic condition. However, all these studies above are for conventional internal combustion engine (ICE) trucks. In this work, we study the problem for electric trucks with a new objective of minimizing the carbon footprint, imposing unique challenges to the problem due to the coupling nature of the State of Charge (SoC) constraints.
See Tab.~\ref{tab:related-truck} for a comparison of our work with related studies.

\textbf{Path and Speed Planning for Electric Vehicle} has long been an active research area. Path planning for EVs involves negative edge weights due to the regenerative system. 
Many variants of the RSP problem have been proposed for EVs by minimizing the energy consumption subject to a deadline~\cite{storandtQuick2012, liuConstrained2017,celaEnergy2014b}, or by finding the fastest path subject to the battery capacity limit~\cite{wangContextAware2013, baumShortest2019}. However, all these works only focus on path planning and ignore the design space of speed planning. More recent works consider both path planning and speed planning~\cite{baumModeling2020, baumSpeedConsumption2014, fontanaOptimal2013, hartmannEnergyEfficient2014,strehler2017energy}. 
For example, Baum et al.~\cite{baumModeling2020} seeks to find the fastest path with capacity limit by optimizing both path planning and speed planning for EVs without charging stops by the tradeoff function propagating technique. Strehler et al.~\cite{strehler2017energy} give theoretical insights to the path and speed planning problem for EVs and provide a fully polynomial time-approximation scheme (FPTAS) for the problem. Overall, all mentioned work and their approaches either do not apply to the CFO problem or incur large time complexity in the CFO problem; see also Table~\ref{tab:related-method} and Appendix~\ref{sec:novelty} for a comparison of conceivable approaches.

\textbf{Variants of Vehicle Routing Problem (VRP).} VRP generalizes the NP-hard traveling salesman problem (TSP), where the operator has a fleet of vehicles to fulfill the requirement of a set of customers. The variants of VRP for EV are considered in
electric vehicle routing problem (EVRP)~\cite{schneiderElectric2014,goekeRouting2015,desaulniersExact2016,montoyaElectric2017}. 
Some other variants of VRP consider the environmental impacts of the routing decisions. For example, the green vehicle routing problem~\cite{erdougan2012green} considers VRP for alternative fuel vehicles.
The pollution routing problem~\cite{bektasPollutionRouting2011} considers minimizing the emission for conventional vehicles. 
The electric arc routing problem (eARP)~\cite{fernandez2022arc} considers minimizing the total travel time on an energy-indexed graph.
As a general observation, variants of VRP only consider path planning and are very complex. The existing methods for VRP do not scale to large-scale network. 



\begin{figure}[!tb]
	\centering
	\includegraphics[width=0.9\linewidth]{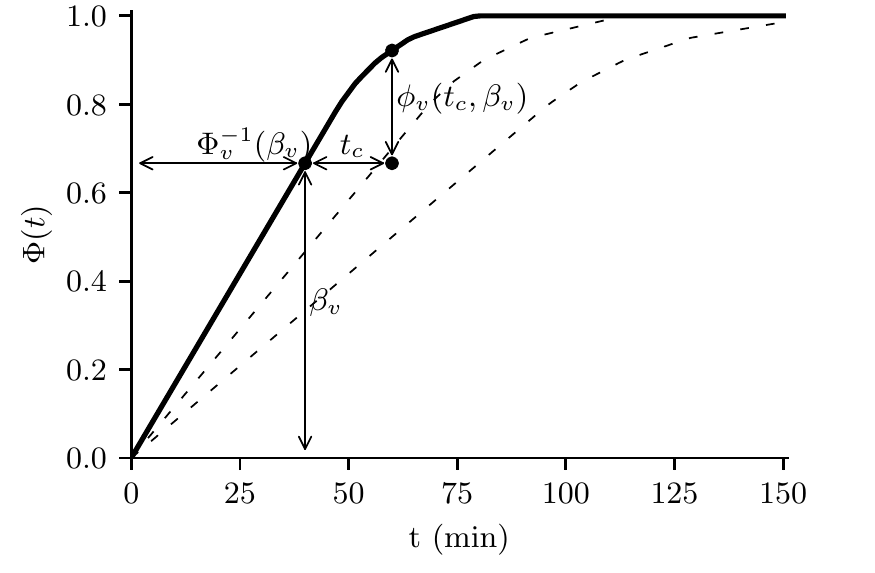}
    \vspace{-4mm}
	\caption{An example of charging function $\Phi$, and the corresponding SoC difference function $\phi$. The initial SoC is $\beta_v=0.67$. After charging for $20$ minutes, 
  the increment of SoC $\phi(t_{c}, \beta_i)=0.25$. The two dashed lines represent charging functions for different charging schemes.}
	\label{fig:charge}
\end{figure}

\section{Problem Formulation}

\subsection{System Model}

\textbf{Transportation Graph and Energy Consumption.}
We model the highway system as a directed graph $ \mathcal{G} = (\mathcal{V}, \mathcal{E}) $ where $\mathcal{V}=\mathcal{V}_r \cup \mathcal{V}_c$ is the set of nodes. 
Here $\mathcal{V}_r$ denotes the nodes for real road segments and $\mathcal{V}_c$ denotes the nodes for the charging stations.
We denote $\mathcal{E} \subset \mathcal{V} \times \mathcal{V}$  the set of edges. 
For each road segment $e \in \mathcal{E}$, we denote its length by $D_{e}$  and  its minimum (resp. maximum) speed limit by $R^{lb}_{e}$ (resp. $R^{ub}_{e}$). We define the minimum and maximum traveling time limit as $t_e^{lb}=D_e/R^{ub}_e$ and $t_e^{ub}=D_e/R^{lb}_e$.


Without loss of generality, we assume homogeneous road condition on a road segment $e \in \mathcal{E}$, e.g., the road grade. Given the weight of an e-truck and the road grade, we model its energy consumption rate over a road segment $e\in \mathcal{E}$ as a convex function of the traveling speed as justified in~\cite{dengEnergyefficient2017}: 
$f_{e}(r_{e}): [R_{e}^{lb}, R_{e}^{ub}] \rightarrow \mathbb{R}$ (unit: kilowatts), where $r_e$ is the traveling speed. Note that for e-trucks, $f_{e}$ can take negative values over (downhill) road segments because of  regenerative breaking~\cite{baumShortest2019}. 
Because of the convexity of $f_{e}(\cdot)$, it suffices to assume that the e-truck travels on the road segment at constant speed~\cite[Lem. 1]{dengEnergyefficient2017} with no optimality loss.\footnote{{We ignore the acceleration/deceleration stage between road segments, as it usually occupies a small fraction of time/energy compared with that under the entire road segment.}} 
We further define an \emph{energy consumption function} for an e-truck to traverse the road segment $e$ as $c_e(t_e) = t_e\cdot f_e(D_e/t_e)$. We note that $c_e(\cdot)$ is a perspective function of $f_e(\cdot)$ with linear transformation. As such $c_e(t_e)$ must be convex in $t_e$~\cite{boyd2004convex}.



\textbf{E-Truck Charging.}
For each charging station node $v \in \mathcal{V}_c$, an e-truck can make two decisions: wait for $t_w \in [ {t}_w^{lb}, t_w^{ub} ]$ amount of time for lower carbon intensity and charge for $t_c \in [0, t_c^{ub}]$ amount of time. Here we use ${t}_w^{lb}$ to denote the time overhead for an e-truck driver to stop and pay bills, etc.
We denote the charging function at $v\in \mathcal{V}_c$ by a concave function $\Phi_{v}(t)$, which represents the charged SoC from zero after we charge for $t$ amount of time.\footnote{{Most e-trucks use lithium-ion batteries, which are often charged with the constant-current constant voltage (CC-CV) scheme~\cite{pelletier2017battery}. In the scheme, a battery is first charged by a constant current, and the SoC increases linearly until  $\sim$80\% of the capacity. After that, a constant voltage is held to avoid overcharging degradation, causing the current to decrease and the SoC to increase concavely.}} Given an initial SoC $\beta_v$  we can compute the increment of SoC after a charging time $t_c$ as 
\begin{equation}
\phi_{v}(t_c, \beta_v) = \Phi_{v}(\Phi^{-1}_{v}(\beta_v) + t_{c}) - \beta_{v}. \label{eq:soc_difference_function}
\end{equation}
Here, $\Phi_v^{-1}$ is the inverse function of $\Phi_v$; see Fig.~\ref{fig:charge} for an illustration.


{There could be multiple charging modes for e-trucks, e.g.,  fast charge and regular charge, which correspond to different charging functions (c.f. Fig.~\ref*{fig:charge}). Here, for simplicity, we assume the e-truck only applies the fast charging scheme at each charging station. We note that our approach can be easily extended to take multiple charging schemes into account. For example, we can model different charging schemes at a charging station by adding multiple charging nodes $i \in \mathcal{V}_c $ at the same location.}

\textbf{Carbon Intensity and Carbon Footprint.} 
For each charging station $v \in \mathcal{V}_c$,
we define a continuous carbon intensity function $\pi_{v}(\tau)$ (unit: kg/kWh)
where $\tau$ is the arrival time instant at node $v$. We assume $\tau_0=0$ at the origin $s$ for ease of the presentation.
As illustrated in Fig.~\ref*{fig:price-distri}, the carbon intensity function $\pi_{v}(\cdot)$ largely varies in locations and time due to the geographical and temporal variation of the renewable generation. 
We note that the carbon-intensity information can be well forecasted at the hourly scale~\cite{majiDACF2022}.
Given the carbon intensity function $\pi_{v}(\cdot)$, if the e-truck with initial SoC $\beta_v$ starts charging at $\tau_v$ at station $v$ and charges $t_c$ time, the induced carbon footprint is

\begin{align}
    F_{v}(\beta_v, t_{c}, \tau_v) = 
    \frac{1}{\eta} \int_{0}^{t_{c}} \pi_{v}(\tau_v + \xi) \frac{\partial_{-} \phi_{v}}{\partial t} (\xi, \beta_v) d\xi. \label{eq:carbon-footprint}
\end{align}
Here, $ 0<\eta \leq 1$ is the charging efficiency of the battery. $\frac{\partial_{-} \phi_{v}}{\partial t}(\xi, \beta)$ is the left partial derivative with respect to $t$, which denotes the charging rate given the charging time $\xi$ and initial SoC $\beta$.
For charging at station $v \in \mathcal{V}_c$, the carbon footprint depends on the initial SoC $\beta_v$, the charging time $t_{c}$ and the arrival time $\tau_v$. It is the integral of the product of the price function and charging rate and thus nonlinear and non-convex (see an example in Fig.~\ref{fig:obj}).

\subsection{Carbon Footprint Optimization Problem}

We study the problem of minimizing the carbon footprint of an e-truck traveling from the origin $s \in \mathcal{V}$ to the destination $d\in \mathcal{V}$, subject to the hard deadline $T$ and the SoC constraints. We consider that the e-truck with initial SoC $\beta_0$ and battery capacity $B$ can take at most $N$ charging stops during the trip. We assume $N$ is a given parameter to our problem. In practice, it is reasonable to limit the number of charging stops in a trip for a good experience. The number is usually small, e.g., 3 or 4 for a 1500-mile e-truck transportation. In our simulation, we set the number of charging stops to be no more than eight. 

\begin{figure}[!tb]
    \centering
    \includegraphics[width=.95\linewidth]{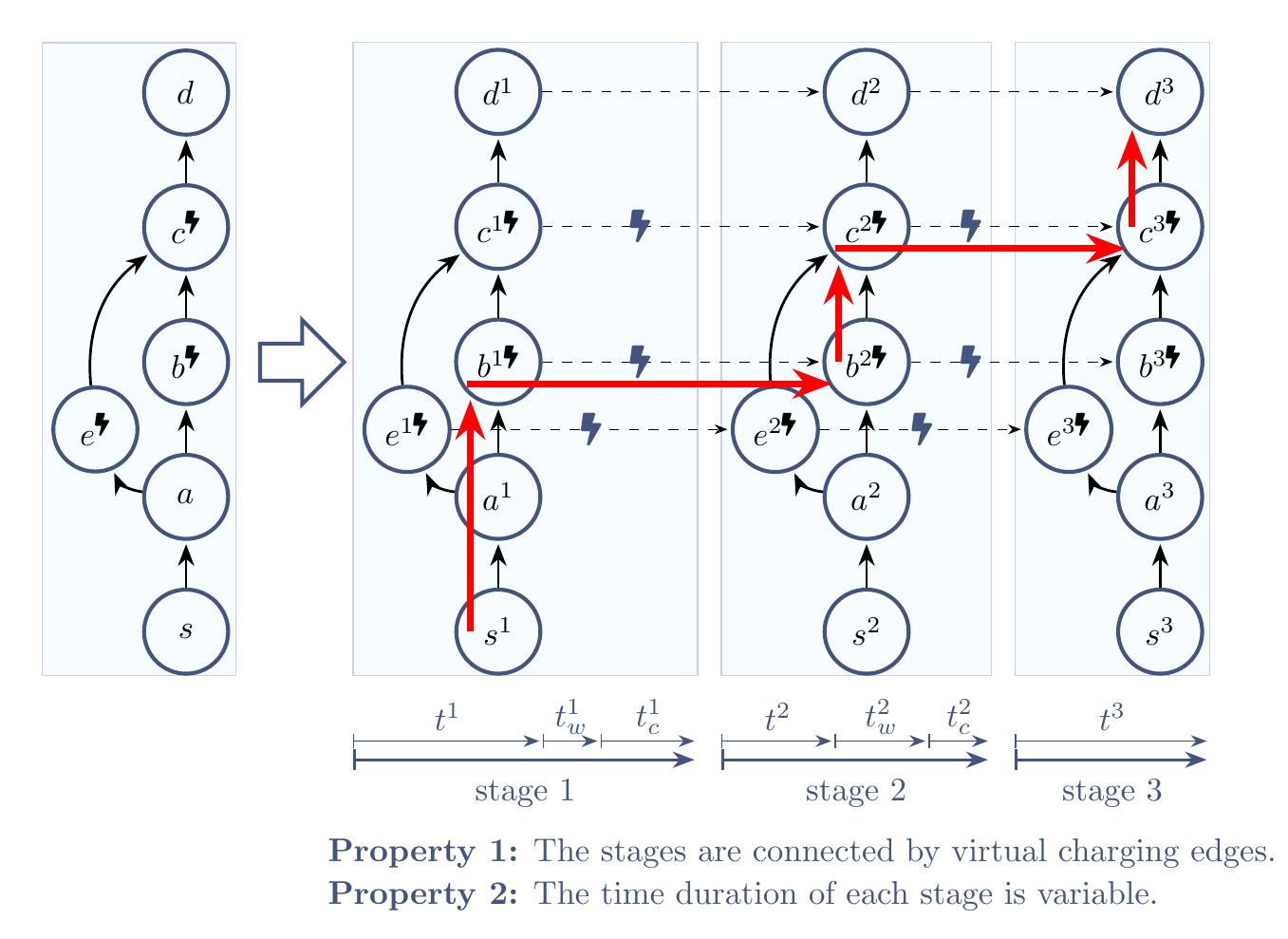}
    \vspace{-3mm}
    \caption{
        An example of stage-expanded graph with $N=2$ charging stops. Here, node $b,c,e$ are charging stations in the original transportation graph. 
    } 
    \label{fig:ex-subpath}
\end{figure}

\textbf{A Stage-Expanded Graph.}
One of the challenges for formulating the CFO problem is the lack of the arrival time information and the arrival SoC information in the original transportation graph $\mathcal{G}$. Those information are essential for objective evaluation in the charging planning. 
To incorporate those information and mitigate this challenge, we consider a \emph{stage-expanded graph}, which extends the original transportation graph to model the charging decisions. 

In particular, we consider a stage-expanded graph $\mathcal{G}_{s}=(\mathcal{V}_s, \mathcal{E}_s)$ with $N+1$ stages. Here $\mathcal{V}_s$ contains $N+1$ copies of original node, i.e.,
$\mathcal{V}_s = \left\{ v^i : v \in \mathcal{V}, i \in \left\{ 1,...,N+1 \right\}  \right\}$. The edge set $\mathcal{E}_s$ keeps the original edges with additional connection between the same charging stations and destinations for two adjacent stages. That is, 
\begin{align*}
    \mathcal{E}_s  = & \left\{ (u^i, v^i) : (u,v) \in \mathcal{E}, i \in \left\{ 1,...,N+1 \right\} \right\}  \\
     & \cup \left\{ (v^i,v^{i+1}) : v \in \mathcal{V}_c \cup \left\{ d \right\}, i \in \left\{ 1,...,N \right\}  \right\}.
\end{align*}
When an e-truck at stage $i$ decides to charge at a charging station $v \in \mathcal{V}_c$, we model it travels through a virtual stage transition edge $(v^i,v^{i+1})$. 
At each stage, the e-truck travels towards the next charging station (or the destination), then wait and charge at the charging station and enter the next stage or stop at the destination.
Note that the e-truck can charge less than $N$ times, arrive the destination at stage $i<N+1$ and then travel across stages through virtual edges between destinations and arrive the final stage.
We provide an illustration in Fig.~\ref{fig:ex-subpath}. 

As we shall see in the following, 
the stage-expanded graph naturally incorporates the information of the arrival time and the arrival SoC, making it easier to model the charging decisions.
We shall also discuss the benefits of such a stage-expanded graph and compare it conceivable alternatives later in Sec.~\ref{sec:challenge} after we formulate the CFO problem.

\textbf{Path and Charging Station Selection.}
{A path in the stage-expanded graph consists of $N+1$ subpaths connecting the source, no more than $N$ charging stations and the destination. We use binary variables to represent the choice of charging stops at each stage and road segments in each subpath. } 

 In particular, for each charging station $v\in \mathcal{V}_c$, we define a binary charging station selection variable $y_v^i \in \left\{ 0,1 \right\}$ and $y_v^i=1$ means we want to choose charging station $v$ as our $i$-th charging stop. For each subpath $i \in \left\{ 1,...,N+1 \right\}$ and road segment $e\in \mathcal{E}$, we define a binary path selection variable $x_e^i \in \left\{ 0,1 \right\}$ and $x_e^i=1$ means we want to travel through road segment $e$ for $i$-th subpath. 
For ease of presentation, we stack the variables 
$\vec{x}^i=\{ x_e^i \}_{e \in \mathcal{E}}$, $\vec{x}=\{ \vec{x}^i \}_{i=1}^{N+1}$ and
$\vec{y}^i=\{ y_e^i \}_{e \in \mathcal{E}}$, $\vec{y}=\{ \vec{y}^i \}_{i=1}^{N+1}$. We define the feasible set of $(\vec{x}, \vec{y})$ as follows:

\begin{subequations}
    \begin{align}
        \mathcal{P} =  \Biggl\{
            & \left( \vec{x}, \vec{y} \right) |\, 
            x_e^i \in \left\{ 0,1 \right\},\ \forall e \in \mathcal{E}, i \in \left\{ 1,...,N+1 \right\} \\
            & y_v^i\in \left\{ 0,1 \right\},\forall  v\in \mathcal{V}, i \in \left\{ 0,...,N+1 \right\} \\
            & \sum_{e \in \text{Out}(v) } x_e^i - \sum_{e \in \text{In}(v) } x_e^i = y_v^{i-1} - y_v^i, \, \forall  v \in \mathcal{V}, i\in\left\{ 1,...,N+1 \right\} \label{eq:pathset:flow} \\
            & \sum_{v\in \mathcal{V}} y_v^i =1, \sum_{v \in \tilde{\mathcal{V}}_c  } y_v^i = 1, \, \forall i \in \left\{ 0,1,...,N+1 \right\} \label{eq:pathset:cs} \\
            & y^0_s = 1, y_d^{N+1} = 1 \label{eq:pathset:bound}
            \Biggr\}.
    \end{align}
\end{subequations}




Here, $\tilde{\mathcal{V}}_c = \mathcal{V}_c \cup \{s,d\} $. Out($v$) and In($v$) are the set of outgoing and incoming edges of node $v$, respectively. Constraint~\eqref{eq:pathset:flow} is the flow conservation requirement for each subpath. Constraint~\eqref{eq:pathset:cs} ensures that only one charging station is selected for each stage. Constraint~\eqref{eq:pathset:bound} defines that boundary case at the origin $s$ and the destination $d$. We remark that  $\mathcal{P}$ also covers path with charging stops less than $N$ as we can choose the same charging stops at adjacent stages.

\textbf{Speed and Charging Planning.}
We then define the decision variables for speed planning and charging planning. For each subpath $i$ and road road segment $e$, we denote the traveling time as $t_e^i$. We also denote by $t_w^{i,v}, t_c^{i,v}$ the waiting and charging time for charging station $v$ as $i$-th charging stop.
We also stack the variables as 
$\vec{t}^i=\left\{ t_e^i \right\}_{e\in \mathcal{E}}$, 
$\vec{t}_w^i=\left\{ t_w^{i,v} \right\}_{v\in \tilde{\mathcal{V}}_c}$, 
$\vec{t}_c^i=\left\{ t_c^{i,v} \right\}_{v\in \tilde{\mathcal{V}}_c}$, 
$\vec{t}=\left\{ \left( \vec{t}^i, \vec{t}^i_w, \vec{t}^i_c \right)  \right\}_{i=0}^{N+1}$. The corresponding feasible set is given by
\begin{subequations}
\begin{align}
    \mathcal{T} = \Biggl\{ 
        \vec{t}\,|\, &t_e^i \in [t_e^{lb}, t_e^{ub}], \forall e\in \mathcal{E}, i \in \left\{ 1,...,N+1 \right\},  \\
        &t_w^{i,v} \in [t_w^{lb}, t_w^{ub}], t_c^{i,v} \in [0, t_c^{ub}], \forall v \in \tilde{\mathcal{V}}_c, i \in \left\{ 1,...,N \right\}, \\
        & t_w^0 = t_c^0 = 0, 
    \Biggr\}.
\end{align}
\end{subequations}

\textbf{Time Constraints.} We denote $\tau_v^i$ the moment of an e-truck entering the charging station $v$ at $i$-th stop. It enables us to allocate the total allowed travel time $T$ to the subpath at each stage. We also stack the variable as $
\vec{\tau}=\left\{ \tau_v^i \right\}_{v\in \tilde{\mathcal{V}}_c, i\in \left\{ 0,...N+1 \right\}}$ and its feasible set is given by

\begin{align}
    \mathcal{T}_\tau = \Biggl\{  \vec{\tau} \,|\,  &\tau_v^i \in [0, T], \forall v \in \tilde{\mathcal{V}}_c, i \in \left\{ 0,...,N+1 \right\} \label{cons:tau:bound}
    \Biggr\}.
\end{align}

We consider the constraint that the traveling time and charging time are within the scheduled time window for $i$-th subpath, $i \in \left\{ 1,...,N+1 \right\}$:

\newcommand{\sumE}{\sum_{e\in \mathcal{E}}}
\newcommand{\sumV}{\sum_{v\in \mathcal{V}_c}}
\begin{equation}
    \begin{aligned}
        \delta_i^\tau(\vec{x}, \vec{y}, \vec{t}, \vec{\tau}) &= \sumE x_e^i t_e^i  + \sumV y_v^{i-1} \left( t_w^{i-1,v} + t_c^{i-1,v} \right) \\ 
    & - \sumV \left( y_v^i \tau_v^i - y_v^{i-1} \tau_v^{i-1} \right) \leq 0.
    \end{aligned}
\end{equation}
Here, the first two terms are the traveling time on $i$-th subpath and the spent time for $(i-1)$-th stop. The third term is the total scheduled time duration between $(i-1)$-th stop and $i$-th stop.

\textbf{SoC Constraints.} It is complicated to keep track of the SoC at each road segment. In practice, we see that the regenerative braking only contributes a small fraction of energy to the e-truck~\cite{su2021Energy}. With such a practical observation, we show that with a small reservation requirement at the SoC when entering a charging stop at a stage, it is sufficient to guarantee the positive SoC in the entire subpath at that stage. 

In particular, we denote by $\beta_v^i$ the SoC entering charging station $v$ at $i$-th stop. We
stack the variable as $
\vec{\beta}=\left\{ \beta_v^i \right\}_{v\in \tilde{\mathcal{V}}_c, i\in \left\{ 0,...N+1 \right\}}$, and its feasible set is given by

\begin{align}
    \mathcal{S}_\alpha = \Biggl\{ \vec{\beta}  \,|\, 
    &\beta_v^i \in [\alpha B, B], \forall v \in \tilde{\mathcal{V}}_c, i \in \left\{ 0,...,N+1 \right\}, \beta_s^0 = \beta_0 \label{cons:betatau:bound}
    \Biggr\}.
\end{align}

Here, $\alpha \in [0,1)$ denotes the ratio of a conservative lower bound for SoC when an e-truck entering a charging station. We consider the initial SoC at the source as $\beta_0$.

We define the SoC constraint that the e-truck does not run out of battery when entering $i$-th stop:
\begin{equation}
    \begin{aligned}
        \delta_i^\beta(\vec{x}, \vec{y}, \vec{t}, \vec{\beta}) &= \sumE x_e^i c_e(t_e^i)  + \sumV y_v^i \beta_v^i \\ 
        & - \sumV y_v^{i-1} \left( \beta_v^{i-1} + \phi_v(\beta_v^{i-1}, t_c^{i-1,v}) \right) \leq 0.
    \end{aligned} \label{cons:soc}
\end{equation}
Here, the first term is the total energy consumption for $i$-th subpath, the second term is the  SoC when the e-truck enters the $i$-th stop, the third term is the SoC when the e-truck leaves the $(i-1)$-th stop. 


We show that a small $\alpha$ is sufficient to ensure SoC feasibility during a subpath, as along as the harvested energy due to regenerative brake is relatively small.
\begin{lemma}\label{lem:alpha:feasible}
    Let an e-truck travels along a subpath with $n$ road segments to a charging stop, along with its energy consumption on each road segment $[c_1,...,c_n]$. If the harvested energy is relatively small and satisfies
    
        \begin{align}
            \frac{1}{2} \sum_{i=1}^{n} \left( |c_i| - c_i \right) & \leq \frac{\alpha}{2(1-\alpha)} \sum_{i=1}^{n} c_i, \label{ineq:recharge} 
        \end{align}
    then the SoC at the charging stop being no less than $\alpha\cdot B$ implies that the e-truck always has a non-negative SoC during this subpath.
\end{lemma}
We leave the proof in Appendix~\ref{sec:lem:alpha:feasible}. Note that the harvested energy (i.e., left hand side of~\eqref{ineq:recharge}) for an e-truck traveling across national highway system is indeed small~\cite{su2021Energy}. Because the maximum permissible grade along the US highway is no larger than $6\%$~\cite{hancock2013policy}, there is little energy that can be harvested.
Meanwhile, our simulation in Sec.~\ref{sec:simulation} shows that $\alpha=0.05$ is sufficient to ensure non-negative SoC for most of the instances in our simulation. {Moreover, it also shows that it incurs a small performance gap by comparing it with a lower bound where $\alpha$ is set as zero. 


\begin{table*}[!tb]
    \centering
    
\newcommand\tikzmark[1]{%
    \tikz[remember picture]  \coordinate[anchor=center] (#1){};%
}
\definecolor{m1color}{rgb}{0.36, 0.54, 0.66}
\definecolor{m2color}{rgb}{1.0, 0.49, 0.0}
\newcommand{\cmark}{\ding{51}}
\newcommand{\xmark}{\ding{55}}%
\resizebox{\linewidth}{!}{
\begin{NiceTabular}{lll| lll}
\toprule
Graph Model in Formulation & Problem & Model Complexity & Algorithm & Complexity & Optimality \\ \midrule
Original graph  & MIP & $ O \left( |\mathcal{V}| + |\mathcal{E}| \right)$ &  \tikzmark{mipAlg}Branch and bound  & \multicellleft{Exponential to \\ $|\mathcal{V}|$ and $|\mathcal{E}|$}  & Optimal \\ \midrule
\multicellleft{Battery-expanded graph} & MILP$^\dagger$ & \tikzmark{milpVar} $ O \left( B / \epsilon \left( |\mathcal{V}| + |\mathcal{E}| \right)\right)$ & Branch and cut  & \multicellleft{Exponential to \\ $|\mathcal{V}|$, $|\mathcal{E}|$, and $B/\epsilon$}& $O(\epsilon)$ to optimal  \tikzmark{milpOptEnd} \\ \midrule
\multicellleft{Time-expanded and \\ battery-expanded graph} &  Shortest Path &  \tikzmark{spVar}  $  \left( (T \cdot B)/\epsilon^2 \right) \left(|\mathcal{V}| + |\mathcal{E}| \right)$ &  Bellman-Ford  & \multicellleft{Polynomial to \\ $|\mathcal{V}|, |\mathcal{E}|, T/\epsilon \text{ and } B/\epsilon$} &  $O(\epsilon)$ to optimal \\ \midrule
\textbf{Our stage-expanded graph} & Generalized RSP$^{\ddagger}$ & \tikzmark{ourVar} $O\big( N \left(|\mathcal{E}| + |\mathcal{V}| \right) \big)$ \tikzmark{ourVarEnd} & \tikzmark{ourAlg}\multicell{Dual subgradient}  & \multicellleft{Polynomial to \\ $|\mathcal{V}|, |\mathcal{E}|, \text{ and } N$} & Posterior bound \tikzmark{ourAlgEnd}\\ \bottomrule 
\end{NiceTabular}
}
\raggedright 
\footnotesize{$^{\dagger}$: Provided that the carbon intensity function is approximated by a piecewise linear function. \quad $^\ddagger$: RSP stands for restricted shortest path.}
\\ \textcolor{m1color}{\textbf{Message 1:}} Compared with battery-expanded graphs, our stage-expanded graph has a low model complexity because th number of stages $N\leq 10$ in practice.
\\ \textcolor{m2color}{\textbf{Message 2:}} Compared with MIP formulations, our formulation reveals an elegant problem structure that leads to an efficient algorithm with favorable performance.

\tikzset{arrow/.style ={draw, line width=1.5pt,dashed}}
\begin{tikzpicture}[remember picture,overlay,->]
    \node[draw=m1color,circle,inner sep=1pt] at ([xshift=1mm, yshift=1mm] ourVarEnd) {\tiny\bf\textcolor{m1color}{1}};
    \draw[m1color,line width=1.2pt] ([xshift=-3mm,yshift=3mm]milpVar.north west) rectangle ([yshift=-1mm,xshift=3mm]ourVarEnd.south);

    \draw[m2color,line width=1.2pt] ([xshift=-6mm,yshift=4mm]mipAlg.north) rectangle ([yshift=-3mm,xshift=-7mm]milpOptEnd.south);
    \node[draw=m2color,circle,inner sep=1pt] at ([xshift=-5mm, yshift=7mm] milpOptEnd) {\tiny\bf\textcolor{m2color}{2}};

    \draw[m2color,line width=1.2pt] ([xshift=-6mm,yshift=6mm]ourAlg.north) rectangle ([yshift=-3mm,xshift=-7mm]ourAlgEnd.south);
    \node[draw=m2color,circle,inner sep=1pt] at ([xshift=-5mm, yshift=1mm] ourAlgEnd) {\tiny\bf\textcolor{m2color}{2}};

\end{tikzpicture}

    \caption{ Comparison of conceivable problem modelings, formulations, and algorithms for the CFO problem. 
    }
    \label{tab:related-method}
\end{table*}

\textbf{Problem Formulation.} To this end, we are ready to formulate the \emph{Carbon Footprint Optimization} (CFO) problem as follows:

\newcommand{\FvObj}{F_v(t_c^{i,v}, t_w^{i,v}, \beta_v^{i}, \tau_v^{i})}

\begin{subequations}
    \begin{align}
      \textsf{CFO}:\;
      \min \;& \sum_{i=1}^{N} \sum_{v \in \mathcal{V}_c}^{} y^{i}_{v} \FvObj \label{obj:cfo:a} \\
        \text{s.t. }\; 
        & \delta_i^\tau(\vec{x}, \vec{y}, \vec{t}, \vec{\tau}) \leq 0,\,  \forall i \in \left\{ 1,...,N+1 \right\}  \label{cons:cfoa:ddl} \\ 
        & \delta_i^\beta(\vec{x}, \vec{y}, \vec{t}, \vec{\beta}) \leq 0 , \, \forall i \in \left\{ 1,...,N+1 \right\} \label{cons:cfoa:soc0} \\
        \text{var.} \; & (\vec{x},\vec{y}) \in \mathcal{P},  \vec{\beta} \in \mathcal{S}_{\alpha},  \vec{\tau}\in \mathcal{T}_\tau , \vec{t} \in \mathcal{T}. \label{cons:cfoa:var}
    \end{align}\label{prob:cfo:a}
\end{subequations}

Here, we aim to minimize the carbon footprint objective~\eqref{obj:cfo:a}. Constraint~\eqref{cons:cfoa:ddl} is the time scheduling constraint for each subpath. Constraints~\eqref{cons:cfoa:soc0} is the SoC constraints for each subpath. We summarize the notations in Tab.~\ref{tab:notation}. 

We give the hardness result of the problem in the follow theorem. 



\begin{theorem}\label{thm:nphard}
    \textsf{CFO} is NP-hard.
\end{theorem}
Please refer to Appendix~\ref{sec:thm:nphard} for the proof. The NP-hardness of the CFO problem comes from the fact that a special case of the CFO problem is PASO problem~\cite{dengEnergyefficient2017} that considers the path planning and speed planning for energy-efficient timely transportation with ICE trucks.

\subsection{Remark on Challenges and Novelty}\label{sec:challenge}

We first remark on four outstanding challenges of the CFO problem, after which we discuss the novelty of our formulation and its significance for tackling those challenges.

\noindent\textbf{Remark on Challenges.}
The first challenge comes from optimal path planning and speed planning under the hard deadline constraint, which already makes the problem NP-hard.  This challenge exists in both the CFO problem and the E2T2 problem for ICE trucks, e.g.,~\cite{dengEnergyefficient2017}. The following three challenges are unique to e-truck and CFO.

The second challenge is the positive battery SoC requirement for each road segment, which adds strong couplings to the problem. One of the conventional ways is to explicitly define variables and constraints on each road segment (c.f., EVRP~\cite{montoyaElectric2017}) and obtain a large Mixed-Integer Programming (MIP). 
However, the resulting MIP is still challenging to solve because general-purpose methods usually do not scale well on MIP.  
We mitigate this challenge in our formulation with a key observation (Lemma~\ref{lem:alpha:feasible}) that it is sufficient to restrict the SoC only when entering each charging stop (and the destination) by a conservative lower bound, i.e.,~\eqref{cons:soc}. As reported in Sec.~\ref{sec:simulation}, such simplification only incurs a small performance loss.

The third challenge comes from the non-convex objective, which is unique to the CFO problem. 
The objective is time-varying and highly non-convex with respect to $\tau$ due to the fluctuating nature of renewable. 
Large-scale non-convex problems are in general challenging to solve. Our problem formulation reveals an elegant problem structure that allows us to decompose the large non-convex problem into a number of significantly less-complex subproblems, which can be efficiently solved by global optimization tools, as we will show in Sec.~\ref{sec:method}.

The last challenge comes from enormous geographical and temporal charging options with diverse carbon intensity. 
The charging station selection configuration space has a combinatorial nature. A straightforward search over locations results in an exponential runtime with respect to $N$. 
To tackle this challenge, our formulation extends the transportation graph to a stage-expanded graph to model charging options. 
This extension significantly reduces the searching space to a polynomial number and allows us to systematically search charging stations with constraint-respect costs via Lagrangian relaxation.

\noindent\textbf{Remark on Novelty.}
Overall, those four challenges make the CFO problem notoriously challenging to solve.
In the following, we remark the novelty of our approach by
comparing our approach with several conceivable alternatives that are applicable to the CFO problem. We summarize the comparison in Tab.~\ref{tab:related-method}. 

The first direction of alternatives (c.f., the first row and the second row in Tab.~\ref{tab:related-method}) resorts to directly formulate the CFO problem as a Mixed Integer Program (MIP) and solve it with general-purpose methods (e.g., branch and bound). However, the resulting MIP is still challenging to solve. The general-purpose methods usually do not scale well to large problems as they do not carefully explore the problem structure. 

The second direction of alternative expands the graph with discretized state $\beta$ and $\tau$ (c.f. the third row in Tab.~\ref{tab:related-method}). Once the expanded graph is constructed, the CFO problem can be transformed to a shortest path problem on the expanded graph. However, the graph is $(T\cdot B)/\epsilon^2$ times of the original graph, making the approach fail to scale to large networks with favorable accuracy.

In contrast, our formulation in~\eqref{prob:cfo:a} enjoys the benefits for the both worlds. It has a low-complexity stage-expanded graph model, which does not involve discretizing the time or SoC. Moreover, the formulation reveals an elegant problem structure for effective algorithm design. 
In the next section, we shall exploit the understanding and design an efficient dual-subgradient method that has good theoretical and empirical performance.
We also provide a detailed comparison in Appendix~\ref{sec:novelty}







\section{An Efficient Algorithm with Performance Guarantee}\label{sec:method}

In this section, we develop an efficient dual-subgradient algorithm for solving the \textsf{CFO} problem with performance guarantee, by exploring the structural insights of the formulation in~\eqref{prob:cfo:a}.

\subsection{The Dual Problem}

We obtain a partially-relaxed dual problem of \textsf{CFO} in~\eqref{prob:cfo:a} by relaxing  the constraints~\eqref{cons:cfoa:ddl} and~\eqref{cons:cfoa:soc0} with dual variable $\lambda_i^\tau$ and $\lambda_i^\beta$, respectively. The Lagrangian function is given by
\begin{align*}
    L(\vec{x}, \vec{y}, \vec{t}, \vec{\beta}, \vec{\tau}, \vec{\lambda}) =&  \sum_{i=1}^{N} \sum_{v\in \mathcal{V}_c}^{} \FvObj  \\
    & + \sum_{i=1}^{N+1} \lambda_i^\tau \delta_i^\tau(\vec{x}, \vec{y}, \vec{t}, \vec{\tau}) 
    +  \sum_{i=1}^{N+1} \lambda_i^\beta \delta_i^\beta(\vec{x}, \vec{y}, \vec{t}, \vec{\beta}),
\end{align*}
where $\vec{\lambda} = (\vec{\lambda}^\beta, \vec{\lambda}^\tau)$. The corresponding dual problem is given by
\begin{align}
    \max_{\vec{\lambda} \geq 0} D ( \vec{\lambda} ) = \max_{\vec{\lambda} \geq 0} 
    \min_{
        \substack{
        (\vec{x},\vec{y}) \in \mathcal{P},  \\ \vec{\beta} \in \mathcal{S}_{\alpha},  \vec{\tau}\in \mathcal{T}_\tau , \vec{t} \in \mathcal{T}
        }} L(\vec{x}, \vec{y}, \vec{t}, \vec{\beta}, \vec{\tau}, \vec{\lambda}).
        \label{prob:dual0}
\end{align}
Given $\vec{\lambda}$, we observe the problem of $D(\vec{\lambda})$ is a two-level problem where the inner level consists of a number of small subproblems and the outer level is a combinatorial problem over the feasible set $\mathcal{P}$. That is, 
\begin{subequations}
    \begin{align}
    D(\vec{\lambda}) & =
         D_1(\vec{\lambda})  
        + \min_{(\vec{x}, \vec{y})  \in \mathcal{P}} 
        \Biggl(  \sum_{i=1}^{N+1} \sumE x_e^i \underbrace{ \min_{t_e^i \in [t_e^{lb}, t_e^{ub}]} g_e^i(\vec{\lambda}, t_e^i) 
        }_{\text{speed planning}} \\
        + &  \sum_{i=1}^{N} \sumV y_v^i 
        \underbrace{
        \min_{\substack{
                t_c^{i,v} \in [0, t_c^{ub}], t_w^{i,v} \in [t_w^{lb}, t_w^{ub}], \\
                \beta_v^i \in [\beta^{lb}, B], \tau_v^i \in [0, T]
            }} 
        h_v^i(\vec{\lambda}, t_c^{i,v}, t_w^{i,v}, \beta_{v}^{i}, \tau_{v}^{i})  
        }_{\text{charging planning}}
        \Biggr) ,
    \end{align} \label{prob:dual}
\end{subequations}
where $D_1(\vec{\lambda}) =  - \lambda_{N+1}^\tau T - \beta_0 \lambda^\beta_1$ is a constant for any given $\vec{\lambda}$ and
\begin{equation}
    g_e^i(\vec{\lambda}, t_e^i)  =  \lambda_i^\tau t_e^i + \lambda_i^\beta c^e(t_e^i)
    \label{func:g}
\end{equation}
represents the trade-off between the traveling time $t_e^i$ and the energy consumption $c_e(t_e^i)$ for each road segment $e\in \mathcal{E}$. The function
\begin{subequations}
    \begin{align}
        & h_v^i(\vec{\lambda}, t_c^{i,v}, t_w^{i,v}, \beta_v^i, \tau_v^i) 
         =  \FvObj + \lambda_{i+1}^\tau \left( t_w^{i,v} + t_c^{i,v} \right) \\
        &\qquad\quad  +  \left( \lambda_{i+1}^\tau - \lambda_{i}^\tau \right) \tau_v^i 
          + \left( \lambda_i^\beta - \lambda_{i+1}^\beta \right) \beta_v^i  
         -\lambda^\beta_{i+1} \phi_v\left(\beta_v^i, t_c^{i,v}\right)
    \end{align} \label{func:h}
\end{subequations}
represents the trade-off between the objective and the constraints. Next, we define 
\begin{subequations}
    \begin{align}
        w_e^i(\vec{\lambda}) &= \min_{t_e^i \in [t_e^{lb}, t_e^{ub}]} g_e^i(\vec{\lambda}, t_e^i), \label{prob:w} \\
        \sigma_v^i(\vec{\lambda}) &=  \min_{\substack{
                    t_c^{i,v} \in [0, t_c^{ub}], t_w^{i,v} \in [t_w^{lb}, t_w^{ub}], \\
                    \beta_v^i \in [\alpha B, B], \tau_v^i \in [0, T]
                }} 
            h_v^i(\vec{\lambda}, t_c^{i,v}, t_w^{i,v}, \beta_{v}^{i}, \tau_{v}^{i}). \label{prob:sigma}
    \end{align}
\end{subequations}
Here, $w_e^i(\vec{\lambda})$ can be computed by solving a single-variable convex optimization problem because $g_e^i(\vec{\lambda}, t_e^i)$ is convex with respect to $t_e^i$. The computation of $\sigma_v^i(\vec{\lambda})$ requires more effort, where we need to solve a 4-variable non-convex optimization problem over box constraints. Fortunately, the dimension of the subproblem~\eqref{prob:sigma} is fixed. Therefore, we can obtain $\sigma(\vec{\lambda})$ in $O( M^4/\epsilon_1^4 )$ with a Branch and Bound (BnB) scheme~\cite{boyd2007branch}. Here $M=\max\{ t_c^{ub}, t_w^{ub}, B, T\}$ is the diameter of the box constraint in subproblem~\eqref{prob:sigma} and $\epsilon_1$ is the accuracy.

\noindent \textbf{Solving the Outer Level Problem.}
After solving the inner level subproblems~\eqref{prob:w} and \eqref{prob:sigma}, the dual function can be written as follows:
\begin{align}
    \hspace{-1mm} D(\vec{\lambda}) = & D_1(\vec{\lambda}) +  \min_{(\vec{x},\vec{y}) \in \mathcal{P}} \sum_{i=1}^{N+1} \sum_{e\in \mathcal{E}}^{} w_e^i(\vec{\lambda}) x_e^i + \sum_{i=1}^{N} \sum_{v \in \mathcal{V}_c}^{} \sigma_v^i(\vec{\lambda}) y^i_v. \label{prob:ilp}
\end{align}
The problem in~\eqref{prob:ilp}, at first glance, is an integer linear program (ILP) which is in general challenging to solve.
However, after exploring the problem structure,
we can see that each charging station $v$ is associated with a cost $\sigma_v^i$ if it is the $i$-th charging stop, and each road segment $e$ is with a cost $w_e^i$ if it is in $i$-th subpath. Our goal is to find $N$ charging stops and $N+1$ subpaths with the minimum total cost. We observe that given any charging station selection, optimal subpaths can be determined by the shortest paths. Correspondingly, we can construct an extended graph consisting of $N$ layers of charging stations with cost $\sigma_v^i$ and edges between any two nodes in the consecutive layers with weight determined by the shortest paths. We can then find the optimal charge stop selection by the shortest path on the extended charging station graph. This idea leads to a novel two-level shortest path algorithm. 

\begin{figure}[!tb]
    \resizebox{.7\linewidth}{!}{
    \begin{tikzpicture}[->,>=stealth',shorten >=1pt,auto,node distance=3cm,
    thick,main node/.style={circle,draw,font=\sffamily\Large\bfseries}]

    \node[main node] (s) {$s$};
    \node[main node] (u1) [above right of=s] {$u^1$} ;
    \node[main node] (u2) [right of=u1] {$u^2$};
    \node[main node] (v1) [below right of=s] {$v^1$} ;
    \node[main node] (v2) [right of=v1] {$v^2$};
    \node[main node] (d) [above right of=v2] {$d$};



    \path[every node/.style={font=\sffamily\small}]
    (s) edge  node [above, rotate=45]  {\footnotesize{$\text{SP}^1(s,u)$}} (u1) 
    (s) edge  node [below, rotate=-45] {\footnotesize{$\text{SP}^1(s,v)$}} (v1) 
    (u1) edge  node [anchor=center, above left, midway, rotate=-55] {\footnotesize{$\sigma^1_u + \text{SP}^2(u,v)$}} (v2) 
    (u1) edge  node [anchor=center, above, midway] {\footnotesize{$\sigma^1_u$}} (u2) 
    (v1) edge  node [anchor=center, above, midway] {\footnotesize{$\sigma^1_v$}} (v2) 
    (v1) edge  node [anchor=center, above left, midway, rotate=55]  {\footnotesize{$\sigma^1_v + \text{SP}^2(v,u)$}} (u2) 
    (u2) edge  node [anchor=center, above, midway, rotate=-45]  {\footnotesize{$\sigma^2_u + \text{SP}^3(u,d)$}} (d) 
    (v2) edge  node [anchor=center, above, midway, rotate=45]  {\footnotesize{$\sigma^2_v + \text{SP}^3(v,d)$}} (d) 
    ;
                
\end{tikzpicture} 
    }
    \vspace*{-1mm}
    \caption{An example of the extended graph with $N=2$. In the original graph, we have two charging stations $u$ and $v$.}
    \vspace*{-1mm}
\end{figure}
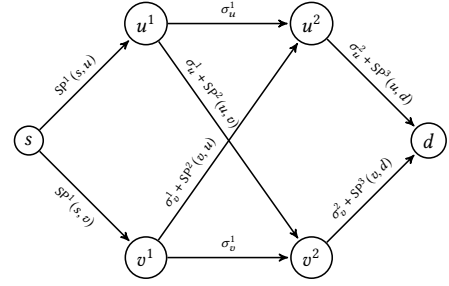
   
In particular, we construct an extended graph $\mathcal{G}_{ex}=(\mathcal{V}_{ex}, \mathcal{E}_{ex})$, where the node set $\mathcal{V}_{ex}$ contains $s,d$ and $N$ copies of charging stations, i.e., 
\begin{align}
    \mathcal{V}_{ex} = \left\{ s,d \right\} \cup \left\{ v^i \, |v \in \mathcal{V}_c, i \in \left\{ 1,...,N \right\} \,   \right\}.
\end{align}
For each $u,v \in \mathcal{V}_c$ and $i\in \left\{ 1,...,N \right\}$, we construct an edge between $u^{i}$ and $v^{i+1}$ with cost 
$ \sigma_{u}^{i}(\vec{\lambda}) + SP^{i+1}(u, v)$
where $SP^{i+1}(u,v)$ is the cost of shortest path between $(u,v)$ on the original graph $\mathcal{G}$ with weight $w_e^{i+1}(\vec{\lambda})$. Similarly, we construct edges $(s,v^1), \forall v \in \mathcal{V}_c$ with weights $SP^1(s,v)$ and construct edges $(u^N, d), \forall v \in \mathcal{V}_c$ with weights $\sigma_v^{N} + SP^{N+1}(u^N, d)$.
By construction, we have the following lemma:

\begin{lemma}
    A shortest path on $\mathcal{G}_{ex}$ corresponds to an optimal solution to ILP problem in~\eqref{prob:ilp}.
\end{lemma}

Hence, we can solve the ILP problem~\eqref{prob:ilp} by applying a shortest path algorithm to the extended graph $\mathcal{G}_{ex}$ with low time complexity. 
Overall, given any $\vec{\lambda}$, we can compute the value of $D(\vec{\lambda})$ by first solving a number of subproblems~\eqref{prob:w} and~\eqref{prob:sigma} and then solving the outer problem~\eqref{prob:ilp} with shortest algorithms on an extended graph.


\subsection{A Dual Subgradient Algorithm}

We use the dual subgradient method to iteratively update the dual variable $\vec{\lambda}$. In particular, we denote by $\vec{\lambda}[k]$ the dual variable at $k$-th iteration. Then we update the dual variable as follows:
\begin{subequations}
    \begin{align}
        \lambda_i^\beta[k+1] &= \Biggl( \lambda_i^\beta[k] + \theta_k 
         \underbrace{ 
            \delta_i^\beta \left(\vec{x}^*[k], \vec{y}^*[k], \vec{t}^*[k], \vec{\beta}^*[k]\right)
        }_{\delta_i^\beta[k]}
          \Biggr)_+,  \\
        \lambda_i^\tau[k+1] &=\Biggl(\lambda_i^\tau[k] + \theta_k 
        \underbrace{
            \delta_i^\tau \left(\vec{x}^*[k], \vec{y}^*[k], \vec{t}^*[k], \vec{\tau}^*[k]\right)  }_{\delta_i^\tau[k]}
        \Biggr)_+.
    \end{align}\label{eq:dualupdate}
\end{subequations}
Here, $\vec{x}^*[k], \vec{y}^*[k], \vec{t}^*[k], \vec{\beta}^*[k], \vec{\tau}^*[k]$ are the solutions of the computed from the subproblems~\eqref{prob:ilp},~\eqref{prob:w}, and~\eqref{prob:sigma} with $\vec{\lambda}[k]$. 
The values $\delta_i^\beta[k]$ and $\delta_i^\tau[k]$ are thus the corresponding subgradients of $D(\vec{\lambda})$ for each components.
The notation $(a)_+$ denotes  $\max\{0,a\}$. The step size at $k$-th iteration is denoted as $\theta_k$.

\begin{algorithm}[tb]
    \caption{A Dual Subgradient Method}\label{alg:dualsub}
    \begin{algorithmic}[1]
        \State \textbf{Initialization:} $\text{sol} \gets \text{NULL}$, $\vec{\lambda}[1] \gets 0$
        \For{$k \gets 1$ to $K$}
            \State Compute $\text{sol}^+$ with $\vec{\lambda}[k]$ by solving problem~\eqref{prob:dual} \label{alg:dualsub:primal}
            \If{$\delta_i^\beta[k] = 0$ and $\delta_i^\tau[k] = 0, \forall i\in \left\{ 1,...,N+1 \right\}$}
                \State \Return sol $\gets \text{sol}^+$  \label{alg:dualsub:optimal}
            \EndIf
            \If{$\delta_i^\beta[k] \leq 0$ and $\delta_i^\tau[k] \leq 0, \forall i\in \left\{ 1,...,N+1 \right\}$ \textbf{and} $\text{sol}^+$ has less objective }
                \State sol $\gets \text{sol}^+$ 
            \EndIf
            \State Compute $\vec{\lambda}[k+1]$ according to~\eqref{eq:dualupdate}.
        \EndFor
        \State \Return sol
    \end{algorithmic}
\end{algorithm}

The overall algorithm is summarized in Algorithm~\ref{alg:dualsub}.
From a high-level perspective, the dual variables $\vec{\lambda}$ can be interpreted as the price for the scheduled deadline and SoC. 
For example, a higher $\lambda^\beta_i$ gives more weight on the constraint $\delta_i^\beta$ and has multiple effects. First, it gives more weight for the energy cost in ~\eqref{func:g}, thus resulting in a larger traveling time in the edge subproblem~\eqref{prob:w}. Meanwhile, a higher $\lambda^\beta_i$ results in a lower scheduled SoC $\beta_v^i$ at $i$-th stop but a higher $\beta_v^{i-1}$ and more charging time $t_c^{i,v}$ at $(i-1)$-th stop. It means that a higher $\lambda^\beta_i$ results in a larger difference of scheduled $SoC$ between charging stops $(i-1)$ and $i$. A similar interpretation applies to the dual variable $\lambda^\tau_i$. 
Overall, Alg.~\ref{alg:dualsub} seeks to find a solution that balances the objective and the constraint cost and adjusts the dual variable via the subgradient direction.


\subsection{Performance Analysis}

In this subsection, as one of the key contributions, we show that (i) our algorithm always converges with a rate of $O(1/\sqrt{K})$ where $K$ is the number of iteration, (ii) each iteration only incurs polynomial-time complexity in graph size, and (iii) it generates optimal results if a condition is met and solutions with bounded loss otherwise.

\noindent \textbf{Convergence Rate.}
Similar to standard subgradient methods~\cite{boyd2003subgradient}, algorithm~\ref{alg:dualsub} also converges to dual optimal value at a rate of $O(1/\sqrt{K})$.

\begin{theorem}\label{thm:converge}
    Let $D^*$ be the optimal dual objective and let $\overline{D}_K$ be the maximum dual value over $K$ iterations in algorithm~\ref{alg:dualsub}. With constant step size $\theta_k = \frac{1}{\sqrt{K}}$, we have,  for some constant $C$,
    \begin{align}
        D^* - \overline{D}_K \leq \frac{C}{\sqrt{K}}.
    \end{align}
   
\end{theorem}
The proof is in Appendix~\ref{sec:thm:converge}. It shows that a constant step size $\theta_k=1/\sqrt{K}$ is sufficient to achieve a convergence rate of $O(1/\sqrt{K})$. However, one may achieve better empirical convergence by adaptively updating the step sizes~\cite{bazaraa1981Choice}. Besides, another improvement can be achieved by modifying the subgradient directions~\cite{camerini1975Improving}.


\noindent \textbf{Time Complexity.}
We summarize the time complexity analysis of Alg.~\ref{alg:dualsub} in the following Proposition.
\begin{proposition}\label{prop:complexity}
    The time complexity per iteration of Alg.~\ref{alg:dualsub} is 
    \begin{align}
        O  \left( \left( N+1 \right) |\mathcal{V}_c| |\mathcal{V}| |\mathcal{E}| + N |\mathcal{V}_c| \frac{M^4}{\epsilon_1^4} \right) \label{eq:complexity}
    \end{align},
    where $M = \left\{ t_c^{ub}, t_w^{ub}, B, T \right\}$ and $\epsilon_1$ is the required accuracy for solving the subproblems~\eqref{prob:sigma}
\end{proposition}
The proof is in Appendix~\ref{sec:prop:complexity}. 
Proposition~\ref{prop:complexity} states that at each iteration , Alg.~\ref{alg:dualsub} incurs only polynomial time complexity in graph size, albeit it requires solving the ILP subproblem in~\eqref{prob:sigma} optimally.
Meanwhile, Theorem~\ref{thm:converge} guarantees that Alg.~\ref{alg:dualsub} converges to a solution with accuracy $\epsilon_0$ in $O(1/\epsilon_0^2)$ iterations. We can derive the total time complexity of Alg.~\ref{alg:dualsub} by combining Proposition~\ref{prop:complexity} and Theorem~\ref{thm:converge}.



\noindent \textbf{Optimality Gap.}
Note that in the dual subgradient method, a convergence to dual optimal does not imply a convergence to the primal optimal as there might be a duality gap.
To this end, we provide the following bound on the optimality gap.
\begin{theorem}\label{thm:optimal}
    Let OPT be the optimal objective to the problem in~\eqref{prob:cfo:a}.
    If Alg.~\ref{alg:dualsub} produces a feasible solution updated at iteration $k$ with objective ALG, then the optimality gap is bounded by
    \begin{align}
        ALG - OPT \leq  - \sum_{i=1}^{N+1} \left( \lambda_i^\beta[k] \delta_i^\beta[k] + \lambda_i^\tau[k] \delta_i^\tau[k] \right). \label{ineq:postbound}
    \end{align}
\end{theorem}

The proof is in Appendix~\ref{sec:thm:optimal}. The posterior bound~\eqref{ineq:postbound} simply comes from the weak duality.
Meanwhile, we can compute it at each iteration in Alg.~\ref{alg:dualsub} and use it for early termination upon certain accuracy.
Moreover, Theorem~\ref{thm:optimal} also provides an optimality condition of the produced solution by the following corollary.
\begin{corollary}
    If Alg.~\ref{alg:dualsub} returns a feasible solution in line~\ref{alg:dualsub:optimal}, then the solution is optimal.
\end{corollary}

Note that Theorem~\ref{thm:optimal} requires Alg.\ref{alg:dualsub} to produce a feasible solution, which is not guaranteed in dual-based methods. In the case where Alg.\ref{alg:dualsub} does not find a feasible solution, we re-optimize the speed planning and charging planning with the given path in the last iteration of Alg.~\ref{alg:dualsub} and try to recover a feasible primal solution. As seen in Sec.~\ref{sec:simulation}, this recovery scheme always finds a feasible solution.



%

\section{Numerical Experiments}\label{sec:simulation}

\newcommand{\algcarbon}{\textsf{CARBON}}
\newcommand{\algenergy}{\textsf{ENERGY}}
\newcommand{\algfast}{\textsf{FAST}}
\newcommand{\algfastspeed}{\textsf{FAST-S}}
\newcommand{\algfastspeedcharge}{\textsf{FAST-SC}}





In this section, we evaluate the performance of our algorithm by simulations using real-world traces over the US highway network. 
Our objectives are 
(i) to study the behavior of the carbon-optimized solutions.
(ii) to study the performance of the proposed approach as compare to  baselines and the impact of path planning, speed planning, charging planning, respectively; 
(iii) to study the impact of conservative ratio $\alpha$; 
(iv) to compare the carbon footprint between electric trucks and ICE trucks, and study the environmental benefit of electric trucks. 

\subsection{Experiment Setup}

\begin{figure}
    \includegraphics[width=0.8\linewidth]{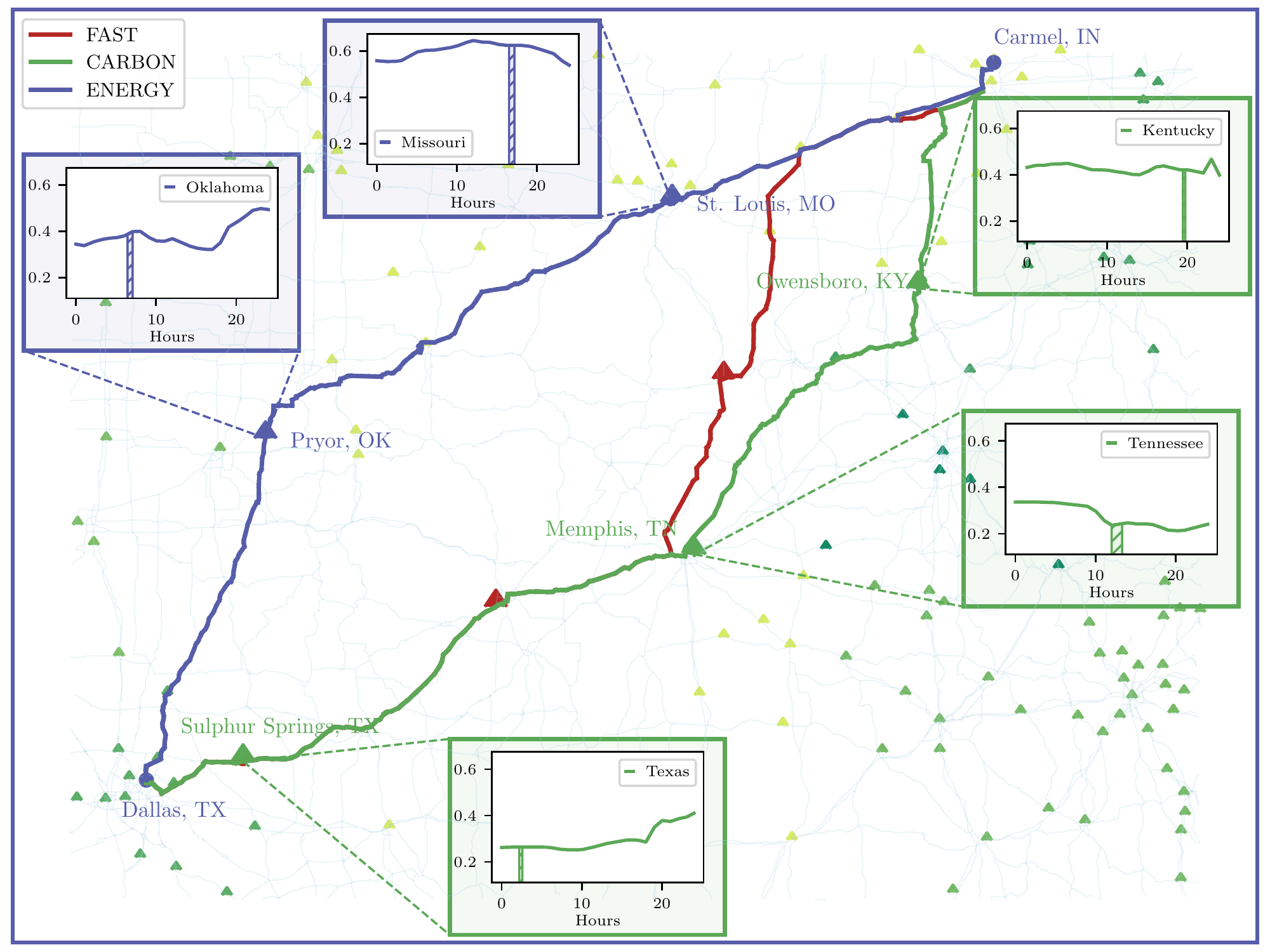}
    \centering
    \resizebox{.8\linewidth}{!}{
    \begin{tabular}{lcccc}
        \toprule
        & \multicell{Carbon \\ footprint (kg)} & \multicell{Energy \\ (kWh)} & \multicell{Distance \\ (miles)} &  \multicell{Time \\ (hours)} \\
        \midrule
        \algfast & 1022.0 & 1633.7  & 919.9& 19.3 \\
        \algenergy & 637.3 & 1239.7 & 879.8& 24.0\\
        \algcarbon & 413.6 & 1487.5 & 946.9& 24.0\\
       \bottomrule
    \end{tabular}
    }
    \caption{The illustration of carbon-optimized and energy-optimized solutions from Dallas, TX to Carmel, IN. 
    In the figure, the shaded area represents the period of charging process.
    } \label{fig:ons-instance}
\end{figure}

\emph{Transportation network.} We collect the highway network data from the Map-based Educational Tools for Algorithm Learning (METAL) project~\cite{metalproject}. 
The constructed network consists of $84,504$ nodes and $178,238$ directed edges. We collect the real-time road speed data from HERE map~\cite{heremap}.
The grade of each road segment is derived from the elevations of its end nodes provided by the Shuttle Radar Topography Mission (SRTM)~\cite{farr2007shuttle} project. 
We also pre-process the road network and merge the adjacent road segments with similar grades with tolerance $0.5\%$ to reduce the computational complexity.
%

\emph{Energy consumption model.}
We use the simulator \emph{FASTSim}~\cite{gonderFuture2018} to collect the energy consumption data with driving speed (from $10$ mph to $70$ mph with a step of $0.2$ mph) at different grades (from $-6\%$ to $6\%$ with a step of $0.25\%$)\footnote{We note that the maximum grade allowed for the inter-state highway is a 6-meter difference in elevation every 100 meters road length in the United States~\cite{hancock2013policy}.}. 
Then we fit the energy consumption data with cubic polynomial functions. 
We use the parameters of Tesla Semi~\cite{teslaSemi} for the electric truck model with the battery capacity $B=1,000$ kWh and the total weight of $36$ tons. We assume the e-truck is fully charged at the origin, i.e., $\beta_0 = B$.

\emph{Origin-destination pair.} We collect origin-destination pairs from the Freight Analysis Framework (FAF)~\cite{hwang2016freight}. We select \numsim origin-destination pairs with distances longer than 800 miles. Those pairs represent 950 billion dollars of freight by trucks in 2017. 
Tab.~\ref{tab:src-des} illustrates a subset of the selected origin-destination pairs in the US starting from Los Angeles.

\begin{table}[tb]
    \centering
    \resizebox{\linewidth}{!}{
    \begin{tabular}{llcc}
        \toprule
        Origin & Destination & \multicell{Distance \\ (miles)} & \multicell{Value \\ (billion USD)}  \\
        \midrule
        Los Angeles CA & Columbus OH & 1977 & 17.725 \\
        Los Angeles CA & Dallas-Fort Worth TX  & 1240 & 12.247 \\
        Los Angeles CA & Chicago IL & 1745 & 11.293 \\
        Los Angeles CA & Nashville TN & 1780 & 10.718 \\
        Los Angeles CA & Houston TX & 1373 & 7.837 \\
       \bottomrule
    \end{tabular}
    }
    \vspace{1mm}
    \caption{Five popular origin-destination pairs from Los Angeles.}\label{tab:src-des}
    \vspace{-10mm}
\end{table}

\emph{Charging station data.}
We collect the location data of charging stations from the OpenStreetMap~(OSM)~\cite{OpenStreetMap} and add each location to the nearest road segment in the constructed graph.
The collected data yields $2,555$ charging stations.
We approximate the charging functions with a piecewise linear function using six break points at $0\%, 80\%, 85\%, 90\%, 95\%, 100\%$ SoC~\cite{baumShortest2019}. 
The resulting charging function can charge the studied electric truck from $0\%$ to $80\%$ in $48$ minutes. See~Fig.\ref{fig:charge} for an illustration. 
We collect the electricity generation data from Energy Information Administration~(EIA)~\cite{eia2021monthly} and obtain a piecewise linear carbon intensity function $\pi(\cdot)$ for each charging station~(See Fig.~\ref*{fig:price-distri}).

\begin{figure*}[tb]
    \includegraphics[width=.95\linewidth]{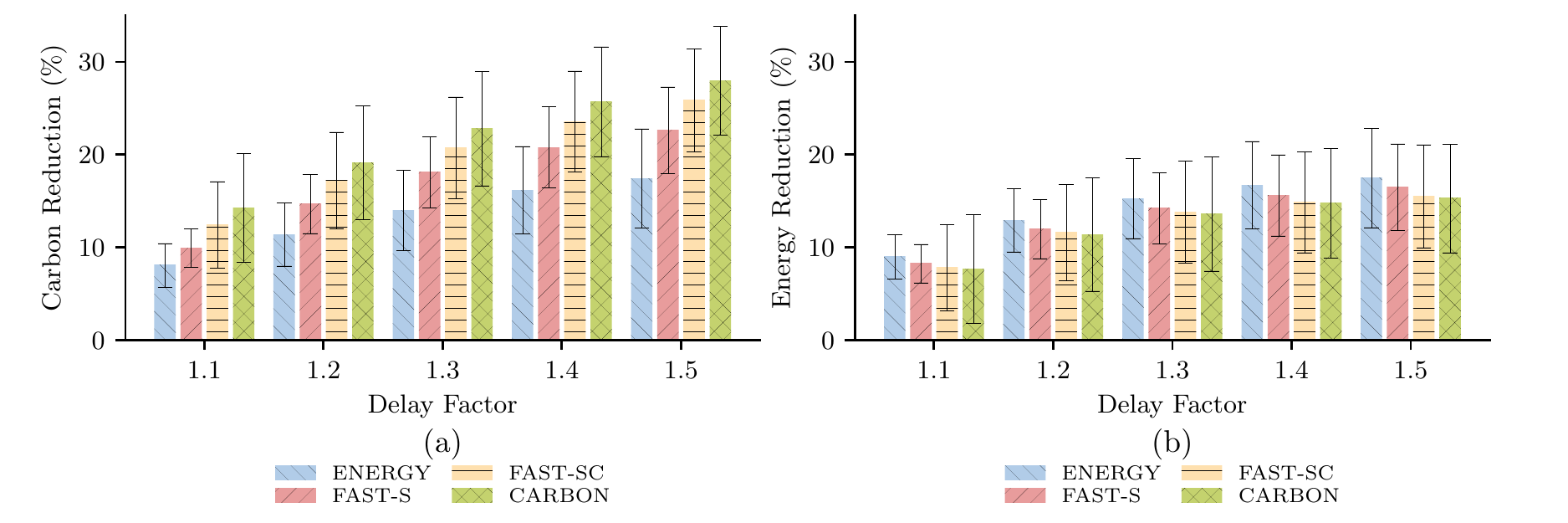}
    \vspace{-5mm}
    \caption{Simulation results for \numsim instances. The carbon (energy) reduction is also computed with respect to \algfast. 
    } \label{fig:baseline}
\end{figure*}


In the simulation, we implement and compare the following conceivable alternatives.

\trianglebullet\ \algcarbon: our carbon footprint minimizing approach. 

\trianglebullet\ \algenergy: our approach with the objective of minimizing energy consumption by setting the carbon intensity function $\pi(\cdot) \equiv 1$.

\trianglebullet\ \algfast: our approach with the objective of minimizing total travel and charging time by changing the objective and following the similar dualization procedure in Sec.~\ref{sec:method}.

\trianglebullet\ \algfastspeed: the same path of FAST with the speed planning.

\trianglebullet\ \algfastspeedcharge: the same path of FAST with the speed planning and charge planning.

\subsection{Behavior of Carbon-Optimized Solution}
We first study a transportation task from Dallas, TX to Carmel, IN. We compare the solution of \algfast, \algcarbon~and \algenergy~and present them in Fig.~\ref{fig:ons-instance}. We observe that \algcarbon~tends to "follow the sun and go with the wind" to charge the electricity with the low carbon intensity and achieves a 59.6\% carbon reduction as compared to \algfast~and 35.2\% as compared to \algenergy.
To achieve such a large carbon reduction, \algcarbon~charge a small amount of energy at the first and the third charging station to ensure the sufficient SoC on road while charge more at the second charging station for lower carbon footprint.
Meanwhile, the detour to the greener charging station leads to a slightly longer total distance and slightly more energy consumption as compared to the energy-optimized solution.
We find similar observations in the next subsection where we conduct the simulation over all 500 origin-destinations.



\subsection{Comparison with Baselines}

We then conduct simulations on \numsim~ origin~destination pairs. For each pair, we denote by $T_f$ the total travel time of \algfast~ and set deadlines from $1.1T_f$ to $1.5T_f$. We call the ratio between the deadline and $T_f$ the \emph{delay factor}.
The start time is set to 4 PM, Jan 1st, 2022. We set the number of charging stops $N$ of \algcarbon~and \algenergy~the same as \algfast. 
{In these typical settings, the number of charging stops is no more than eight.} 
We set the ratio of conservative lower bound $\alpha=0.05$. 
\algcarbon~produces feasible solutions for $96\%$ of all instances. 
For infeasible solutions, we increase $\alpha$ until it produces a feasible solution. A maximum value of $\alpha=0.12$ guarantees \algcarbon~ to produce feasible solutions for all instances.



We show the simulation results in Fig.~\ref{fig:baseline}.
All the results are mean values over \numsim origin-destination pairs and the error bars represent the standard deviations. Fig.~\ref{fig:baseline}(a) and Fig.~\ref{fig:baseline}(b) give the the carbon (resp. energy) reduction of different baselines as compared to \algfast. We observe that as the deadline gets relaxed, the carbon (and energy) reduction of all the baselines increases. 
Moreover, different components of the design space contribute differently into such reduction as we shall discuss in the following.

\noindent\textbf{Benefit of Carbon-Optimized Solution.}
We first compare \algcarbon~and \algenergy.
As shown in Fig.~\ref{fig:baseline}(a), while the carbon footprint of both approaches decrease as the deadline increases, \algcarbon~benefits more from the relaxed deadline and achieves higher carbon reduction.
In particular, , \algcarbon~achieves {19\%} carbon reduction as compared to \algfast~baseline under the typical setting when the delay factor is $1.2$ and up to {28\%} when the delay factor is $1.5$.
Meanwhile, as compared to its energy-efficient counterpart \algenergy, \algcarbon~achieves {9\%} carbon reduction when the delay factor is $1.2$ and {13\%} when the delay factor is $1.5$.
This difference mainly comes from the fact that minimizing energy simply ignores the temporal diversity of the renewable and loses part of the optimization space.
Therefore, minimizing energy consumption does not necessarily lead to a minimum carbon footprint. It is essential to consider the carbon footprint objective when studying the environmental impact of e-trucks.
Meanwhile, we also observe from Fig.~\ref{fig:baseline}(b) that the difference of energy reduction between \algcarbon~and \algenergy~is less than {3\%}.
This means that a carbon-efficient route is also energy-efficient.

\noindent\textbf{Impact of Speed Planning.} We then investigate \algfastspeed~to study the impact of speed planning. When the delay factor is $1.2$. \algfastspeed~saves {15\%} carbon footprint and {12\%} energy as compared to \algfast. As compared to \algenergy, \algfastspeed~saves {3\%} additional carbon footprint. This additional saving comes from the fact that the optimized speed plan re-schedules the travel time of each stage such that the e-truck charges greener electricity at each stage.

\noindent\textbf{Impact of Charging Planning.} We then compare \algfastspeed~and \algfastspeedcharge~to study the impact of charging planning. We observe that when the delay factor is $1.2$, \algfastspeedcharge~saves {2.5\%} additional carbon footprint as compared to \algfastspeed ~because \algfastspeedcharge~takes the advantage of charging planning such that the e-truck can wait for the lower carbon intensity.

\noindent\textbf{Impact of Path Planning.} We then compare \algfastspeedcharge~and \algcarbon~to study the impact of path planning.
We observe that when the delay factor is $1.2$, \algcarbon~saves {2\%} additional carbon footprint as compared to \algfastspeedcharge~because e-truck can select greener charging stations with path planning.

\subsection{Impact of Conservative Ratio $\alpha$}

\begin{figure}[tb]
    \centering
    \includegraphics[width=.7\linewidth]{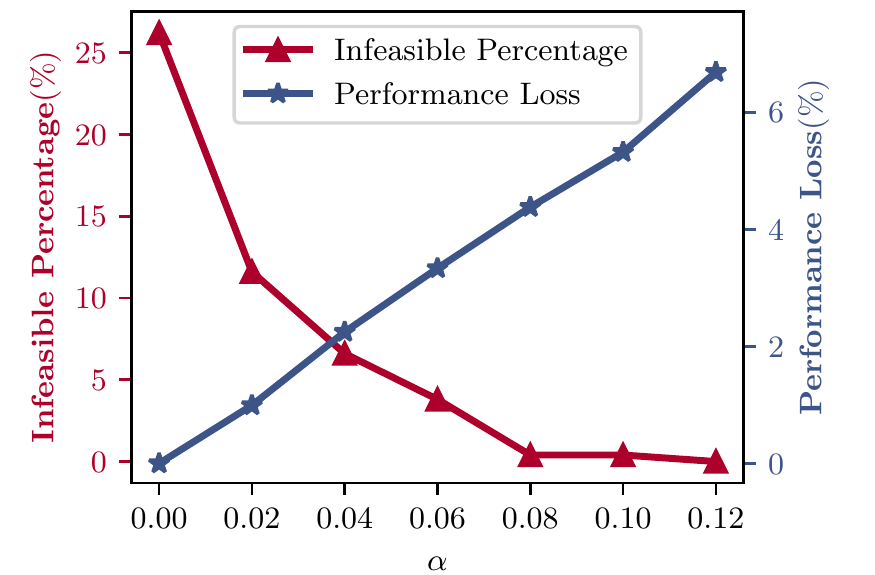}
    \vspace*{-4mm}
    \caption{Impact of conservative ratio $\alpha$. 
    }\label{fig:alpha:performance}
\end{figure}

In this subsection, we study the impact of the conservative ratio $\alpha$. We run \algcarbon~on \numsim origin-destination pairs for $\alpha$ ranging from $0.0$ to $0.12$. 
{For each $\alpha$, we check the SoC at each road segment of the solutions and report the percentage of  solutions violating positive SoC among all \numsim origin-destination pairs. In addition, we notice that \textsf{CFO} with $\alpha=0$ (named $\textsf{CFO}_{0}$) provides a lower bound to the \textsf{CFO} with positive SoC constraints at all road segments but no reservation requirement. We thus define performance loss of \algcarbon~with \textsf{CFO} as its relative error with respect to the \algcarbon~with $\textsf{CFO}_{0}$, i.e.,
$$
\text{Performance Loss} = \frac{ALG_{\alpha}-ALG_{0}}{ALG_{0}},
$$
where $ALG_{\alpha}$ is the objective value of solutions by \algcarbon~for \textsf{CFO}.
The results are presented in Fig.~\ref{fig:alpha:performance}. }
We observe that when $\alpha=0$, there are {26\%} of infeasible solutions while a small increment of $\alpha$ leads to a large feasibility improvement, e.g., $\alpha=0.02$ with {12\%} of infeasible solutions. 
We observe that $\alpha=0.12$ is sufficient to guarantee feasibility for all \numsim origin-destination pairs. Meanwhile, we observe the performance loss introduced by $\alpha$ grows almost linearly with respect to $\alpha$. In a typical selection when $\alpha=0.06$, $96\%$ of solutions are feasible while the performance loss is only {$3.3\%$}. In the most conservative case when $\alpha=0.12$ that guarantees the feasibility of all the solutions, the performance loss is {$6.7\%$}.
Thus, by introducing the conservative ratio $\alpha$, we significantly simplify the coupling SoC constraint while ensuring feasibility with minor performance loss.


\subsection{Comparison with ICE Trucks}

In this subsection, we study the environmental benefits of electric trucks as compared to ICE trucks. 
We use the parameters of Kenworth T800 trailer~\cite{kenworthT800} for the ICE truck model with the same total weight of $36$ tons. 
For the ICE truck, we compute the energy-efficient path plan and speed plan by the dual-based method in~\cite{dengEnergyefficient2017}.
We use the \coo\ factor value of diesel  from EIA~\cite{eiaCO2Factors} to transform the fuel consumption of the ICE truck to its direct carbon emission. The data suggest that if an ICE truck is powered by diesel, it will produce 9.88 kg \coo\ per gallon, or equivalently 0.25 kg \coo\ emission per kWh.
The simulation results are presented in Tab.~\ref{tab:ice:cmp}. We observe that the e-truck saves 55.9\% carbon as compared to the ICE truck under a delay factor $1.2$.

There are several factors that contribute to such carbon reduction. First, e-trucks consume less energy due to the high efficiency of electric motors and the regenerative braking system. 
Second, with the increasing integration of renewable generations for e-truck charging, the carbon footprint of an electric truck reduces. 
Third, our joint optimization of path planning, speed planning and charging planning leads an e-truck to charge clean electricity with low carbon intensity. 
Overall, these factors together lead to a significant reduction of the carbon footprint of the electric truck as compared to the ICE truck.


\begin{table}[tb]
    \centering
    \begin{tabular}{lccccc}
        \toprule
        Delay Factor & 1.1 & 1.2 & 1.3 & 1.4 & 1.5 \\
        \midrule
        \algenergy  & 50.8\% & 51.4\% & 52.4\% & 53.6\% & 54.1\% \\
        \algcarbon  & 54.7\% & 55.9\% & 57.1\% & 58.3\% & 59.2\% \\
        \bottomrule
    \end{tabular}
    \vspace{1mm}
    \caption{Average carbon reduction with respect to ICE truck over \numsim~origin-destination pairs.}\label{tab:ice:cmp}
    \vspace{-6mm}
\end{table}

\section{Discussion}


\noindent\textbf{Battery Degradation.} The modeling and prediction of battery degradation is in general difficult and still an active research area~\cite{hu2020battery}. Among many factors that affects the battery degradation, the most relevant ones for the transportation tasks are charge or discharge rate and the depth of discharge (DoD)~\cite{pelletier2017battery}. 
A high charge rate or a high DoD may accelerate the battery degradation.
To prolong the battery lifetime, we can add a battery degradation cost in the objective. 
We can model different charging mode by adding multiple charging stations at the same location, where a charging mode with the higher rate incur higher degradation cost.
To prevent large DoD, we can also add a degradation cost if the arrival SoC is lower than certain value.


\noindent\textbf{Dynamic Traffic Condition.}
%
In the real world, the traffic conditions can be time-varying. To adapt variable traffic conditions in our formulation, we can incorporate the phase-based traffic model~\cite{xuRide2019} into our stage-expanded graph, where the traffic condition is approximately static at each phase. 
Then we can treat the combination of the stage and the traffic phase as the new stage and apply our approach.

\noindent\textbf{Applications Beyond CFO.}
We remark that our approach has applications beyond this paper. For example, it is also applicable to general routing problem for electric/ICE vehicles with energy consumption and monetary objective. 
Our approach is also applicable to the subproblem for electric vehicle routing problem (EVRP) in a branch-price-and-cut framework~\cite{desaulniersExact2016}.
\section{Conclusion}

In this paper, we present the first study in the literature on the \emph{Carbon Footprint Optimization} (CFO) problem of a heavy-duty e-truck traveling from an origin to a destination across a national highway network subject to a hard deadline, by optimizing path planning, speed planning, and intermediary charging planning.
We show that the problem is NP-hard and under practical settings it is equivalent to finding a generalized RSP on a novel stage-expanded graph, which extends the original transportation graph to model charging options. 
Our problem formulation incurs low model complexity and reveals a problem structure that leads us to an efficient dual-subgradient algorithm with a strong performance guarantee.
Extensive simulation results over the US national highway network show that our solutions reduce up to 28\% carbon footprint compared to the baseline alternatives.
The results also demonstrate the carbon efficiency of e-trucks over ICE trucks.

\begin{acks}
The authors would like to thank the anonymous reviewers for their insightful comments that help further improve this work. 
The work presented in this paper was supported in part by a General Research Fund from Research Grants Council, Hong Kong (Project No. 14207520).
\end{acks}

\bibliography{ref}


\begin{thebibliography}{54}


\ifx \showCODEN    \undefined \def \showCODEN     #1{\unskip}     \fi
\ifx \showDOI      \undefined \def \showDOI       #1{#1}\fi
\ifx \showISBNx    \undefined \def \showISBNx     #1{\unskip}     \fi
\ifx \showISBNxiii \undefined \def \showISBNxiii  #1{\unskip}     \fi
\ifx \showISSN     \undefined \def \showISSN      #1{\unskip}     \fi
\ifx \showLCCN     \undefined \def \showLCCN      #1{\unskip}     \fi
\ifx \shownote     \undefined \def \shownote      #1{#1}          \fi
\ifx \showarticletitle \undefined \def \showarticletitle #1{#1}   \fi
\ifx \showURL      \undefined \def \showURL       {\relax}        \fi
\providecommand\bibfield[2]{#2}
\providecommand\bibinfo[2]{#2}
\providecommand\natexlab[1]{#1}
\providecommand\showeprint[2][]{arXiv:#2}

\bibitem[Administration(2022a)]%
        {eiaCO2Factors}
\bibfield{author}{\bibinfo{person}{U.S. Energy~Information Administration}.}
  \bibinfo{year}{2022}\natexlab{a}.
\newblock \bibinfo{booktitle}{\emph{{Carbon Dioxide Emissions Coefficients}.
  \url{https://www.eia.gov/environment/emissions/co2_vol_mass.php}}}.
\newblock
\newblock
\shownote{Accessed 2022-03-14}.


\bibitem[Administration(2022b)]%
        {energyinformationadministrationeiaRealtime}
\bibfield{author}{\bibinfo{person}{U.S. Energy~Information Administration}.}
  \bibinfo{year}{2022}\natexlab{b}.
\newblock \bibinfo{booktitle}{\emph{Real-Time Operating Grid at
  U.S..~\url{https://www.eia.gov/electricity/gridmonitor/index.php}}}.
\newblock
\newblock
\shownote{Accessed: 2022-03-14}.


\bibitem[Amazon(2022)]%
        {amazon}
\bibfield{author}{\bibinfo{person}{Amazon}.} \bibinfo{year}{2022}\natexlab{}.
\newblock \bibinfo{booktitle}{\emph{Place an Order with Guaranteed Delivery.
  \url{https://www.amazon.com/gp/help/customer/display.html?nodeId=GUDRHL982ZQE9EGL}}}.
\newblock
\newblock
\shownote{Accessed 2022-03-14}.


\bibitem[Ashby(1987)]%
        {ashby1987protecting}
\bibfield{author}{\bibinfo{person}{B~Hunt Ashby}.}
  \bibinfo{year}{1987}\natexlab{}.
\newblock \bibinfo{booktitle}{\emph{Protecting perishable foods during
  transport by truck}}.
\newblock Number 669. \bibinfo{publisher}{US Department of Agriculture, Office
  of Transportation}.
\newblock


\bibitem[Associations(2021)]%
        {truckData2021}
\bibfield{author}{\bibinfo{person}{American~Trucking Associations}.}
  \bibinfo{year}{2021}\natexlab{}.
\newblock \bibinfo{title}{Economics and Industry Data}.
\newblock
  \bibinfo{howpublished}{\url{https://trucking.org/economics-and-industry-data}}.
\newblock
\newblock
\shownote{Accessed: 2022-03-14}.


\bibitem[Baum et~al\mbox{.}(2019)]%
        {baumShortest2019}
\bibfield{author}{\bibinfo{person}{Moritz Baum}, \bibinfo{person}{Julian
  Dibbelt}, \bibinfo{person}{Andreas Gemsa}, \bibinfo{person}{Dorothea Wagner},
  {and} \bibinfo{person}{Tobias Z{\"u}ndorf}.} \bibinfo{year}{2019}\natexlab{}.
\newblock \showarticletitle{Shortest {{Feasible Paths}} with {{Charging Stops}}
  for {{Battery Electric Vehicles}}}.
\newblock \bibinfo{journal}{\emph{Transportation Science}}
  \bibinfo{volume}{53} (\bibinfo{year}{2019}), \bibinfo{pages}{1627--1655}.
\newblock
\showISSN{0041-1655, 1526-5447}


\bibitem[Baum et~al\mbox{.}(2014)]%
        {baumSpeedConsumption2014}
\bibfield{author}{\bibinfo{person}{Moritz Baum}, \bibinfo{person}{Julian
  Dibbelt}, \bibinfo{person}{Lorenz {H{\"u}bschle-Schneider}},
  \bibinfo{person}{Thomas Pajor}, {and} \bibinfo{person}{Dorothea Wagner}.}
  \bibinfo{year}{2014}\natexlab{}.
\newblock \showarticletitle{Speed-{{Consumption Tradeoff}} for {{Electric
  Vehicle Route Planning}}}. In \bibinfo{booktitle}{\emph{{{14Th}} Workshop on
  Algorithmic Approaches for Transportation Modelling, Optimization, and
  Systems}}. \bibinfo{publisher}{{Schloss Dagstuhl - Leibniz-Zentrum fuer
  Informatik GmbH, Wadern/Saarbruecken, Germany}}, \bibinfo{pages}{14 pages}.
\newblock


\bibitem[Baum et~al\mbox{.}(2020)]%
        {baumModeling2020}
\bibfield{author}{\bibinfo{person}{Moritz Baum}, \bibinfo{person}{Julian
  Dibbelt}, \bibinfo{person}{Dorothea Wagner}, {and} \bibinfo{person}{Tobias
  Z{\"u}ndorf}.} \bibinfo{year}{2020}\natexlab{}.
\newblock \showarticletitle{Modeling and {{Engineering Constrained Shortest
  Path Algorithms}} for {{Battery Electric Vehicles}}}.
\newblock \bibinfo{journal}{\emph{Transportation Science}}
  \bibinfo{volume}{54} (\bibinfo{year}{2020}), \bibinfo{pages}{1571--1600}.
\newblock
\showISSN{0041-1655, 1526-5447}


\bibitem[Bazaraa and Sherali(1981)]%
        {bazaraa1981Choice}
\bibfield{author}{\bibinfo{person}{Mokhtar~S Bazaraa} {and}
  \bibinfo{person}{Hanif~D Sherali}.} \bibinfo{year}{1981}\natexlab{}.
\newblock \showarticletitle{On the choice of step size in subgradient
  optimization}.
\newblock \bibinfo{journal}{\emph{European Journal of Operational Research}}
  \bibinfo{volume}{7}, \bibinfo{number}{4} (\bibinfo{year}{1981}),
  \bibinfo{pages}{380--388}.
\newblock


\bibitem[Bekta{\c s} and Laporte(2011)]%
        {bektasPollutionRouting2011}
\bibfield{author}{\bibinfo{person}{Tolga Bekta{\c s}} {and}
  \bibinfo{person}{Gilbert Laporte}.} \bibinfo{year}{2011}\natexlab{}.
\newblock \showarticletitle{The {{Pollution-Routing Problem}}}.
\newblock \bibinfo{journal}{\emph{Transportation Research Part B:
  Methodological}}  \bibinfo{volume}{45} (\bibinfo{year}{2011}),
  \bibinfo{pages}{1232--1250}.
\newblock
\showISSN{01912615}


\bibitem[Bibra et~al\mbox{.}(2022)]%
        {bibra2022global}
\bibfield{author}{\bibinfo{person}{Ekta~Meena Bibra},
  \bibinfo{person}{Elizabeth Connelly}, \bibinfo{person}{Shobhan Dhir},
  \bibinfo{person}{Michael Drtil}, \bibinfo{person}{Pauline Henriot},
  \bibinfo{person}{Inchan Hwang}, \bibinfo{person}{Jean-Baptiste Le~Marois},
  \bibinfo{person}{Sarah McBain}, \bibinfo{person}{Leonardo Paoli}, {and}
  \bibinfo{person}{Jacob Teter}.} \bibinfo{year}{2022}\natexlab{}.
\newblock \showarticletitle{Global EV outlook 2022: Securing supplies for an
  electric future}.
\newblock  (\bibinfo{year}{2022}).
\newblock


\bibitem[Boyd et~al\mbox{.}(2004)]%
        {boyd2004convex}
\bibfield{author}{\bibinfo{person}{Stephen Boyd}, \bibinfo{person}{Stephen~P
  Boyd}, {and} \bibinfo{person}{Lieven Vandenberghe}.}
  \bibinfo{year}{2004}\natexlab{}.
\newblock \bibinfo{booktitle}{\emph{Convex optimization}}.
\newblock \bibinfo{publisher}{Cambridge university press}.
\newblock


\bibitem[Boyd and Mattingley(2007)]%
        {boyd2007branch}
\bibfield{author}{\bibinfo{person}{Stephen Boyd} {and} \bibinfo{person}{Jacob
  Mattingley}.} \bibinfo{year}{2007}\natexlab{}.
\newblock \showarticletitle{Branch and bound methods}.
\newblock \bibinfo{journal}{\emph{Notes for EE364b, Stanford University}}
  \bibinfo{volume}{2006} (\bibinfo{year}{2007}), \bibinfo{pages}{07}.
\newblock


\bibitem[Boyd et~al\mbox{.}(2003)]%
        {boyd2003subgradient}
\bibfield{author}{\bibinfo{person}{Stephen Boyd}, \bibinfo{person}{Lin Xiao},
  {and} \bibinfo{person}{Almir Mutapcic}.} \bibinfo{year}{2003}\natexlab{}.
\newblock \showarticletitle{Subgradient methods}.
\newblock \bibinfo{journal}{\emph{lecture notes of EE392o, Stanford University,
  Autumn Quarter}}  \bibinfo{volume}{2004} (\bibinfo{year}{2003}),
  \bibinfo{pages}{2004--2005}.
\newblock


\bibitem[Camerini et~al\mbox{.}(1975)]%
        {camerini1975Improving}
\bibfield{author}{\bibinfo{person}{P.~M. Camerini}, \bibinfo{person}{L.
  Fratta}, {and} \bibinfo{person}{F. Maffioli}.}
  \bibinfo{year}{1975}\natexlab{}.
\newblock \showarticletitle{On improving relaxation methods by modified
  gradient techniques}.
\newblock In \bibinfo{booktitle}{\emph{Nondifferentiable Optimization}},
  \bibfield{editor}{\bibinfo{person}{M.~L. Balinski} {and}
  \bibinfo{person}{Philip Wolfe}} (Eds.). Vol.~\bibinfo{volume}{3}.
  \bibinfo{publisher}{Springer Berlin Heidelberg}, \bibinfo{pages}{26--34}.
\newblock
\showISBNx{978-3-642-00763-7 978-3-642-00764-4}
\urldef\tempurl%
\url{https://doi.org/10.1007/BFb0120697}
\showDOI{\tempurl}
\newblock
\shownote{Series Title: Mathematical Programming Studies}.


\bibitem[Cela et~al\mbox{.}(2014)]%
        {celaEnergy2014b}
\bibfield{author}{\bibinfo{person}{Arben Cela}, \bibinfo{person}{Tomas Jurik},
  \bibinfo{person}{Redha Hamouche}, \bibinfo{person}{Rene Natowicz},
  \bibinfo{person}{Abdelatif Reama}, \bibinfo{person}{Silviu-Iulian Niculescu},
  {and} \bibinfo{person}{Jerome Julien}.} \bibinfo{year}{2014}\natexlab{}.
\newblock \showarticletitle{Energy {{Optimal Real-Time Navigation System}}}.
\newblock \bibinfo{journal}{\emph{IEEE Intelligent Transportation Systems
  Magazine}}  \bibinfo{volume}{6} (\bibinfo{year}{2014}),
  \bibinfo{pages}{66--79}.
\newblock
\showISSN{1941-1197}


\bibitem[Cormen et~al\mbox{.}(2022)]%
        {cormen2022introduction}
\bibfield{author}{\bibinfo{person}{Thomas~H Cormen}, \bibinfo{person}{Charles~E
  Leiserson}, \bibinfo{person}{Ronald~L Rivest}, {and}
  \bibinfo{person}{Clifford Stein}.} \bibinfo{year}{2022}\natexlab{}.
\newblock \bibinfo{booktitle}{\emph{Introduction to algorithms}}.
\newblock \bibinfo{publisher}{MIT press}.
\newblock


\bibitem[Davis and Boundy(2021)]%
        {davisTransportation2021}
\bibfield{author}{\bibinfo{person}{Stacy Davis} {and}
  \bibinfo{person}{Robert~Gary Boundy}.} \bibinfo{year}{2021}\natexlab{}.
\newblock \bibinfo{booktitle}{\emph{Transportation {Energy} {Data} {Book}:
  {Edition} 39}}.
\newblock \bibinfo{type}{{T}echnical {R}eport}. \bibinfo{institution}{Oak Ridge
  National Lab. (ORNL), Oak Ridge, TN (United States)}.
\newblock
\urldef\tempurl%
\url{https://doi.org/10.2172/1767864}
\showDOI{\tempurl}


\bibitem[Deng et~al\mbox{.}(2017)]%
        {dengEnergyefficient2017}
\bibfield{author}{\bibinfo{person}{Lei Deng}, \bibinfo{person}{Mohammad~H.
  Hajiesmaili}, \bibinfo{person}{Minghua Chen}, {and} \bibinfo{person}{Haibo
  Zeng}.} \bibinfo{year}{2017}\natexlab{}.
\newblock \showarticletitle{Energy-Efficient Timely Transportation of Long-Haul
  Heavy-Duty Trucks}.
\newblock \bibinfo{journal}{\emph{IEEE Transactions on Intelligent
  Transportation Systems}}  \bibinfo{volume}{19} (\bibinfo{year}{2017}),
  \bibinfo{pages}{2099--2113}.
\newblock


\bibitem[Desaulniers et~al\mbox{.}(2016)]%
        {desaulniersExact2016}
\bibfield{author}{\bibinfo{person}{Guy Desaulniers}, \bibinfo{person}{Fausto
  Errico}, \bibinfo{person}{Stefan Irnich}, {and} \bibinfo{person}{Michael
  Schneider}.} \bibinfo{year}{2016}\natexlab{}.
\newblock \showarticletitle{Exact {{Algorithms}} for {{Electric Vehicle-Routing
  Problems}} with {{Time Windows}}}.
\newblock \bibinfo{journal}{\emph{Operations Research}}  \bibinfo{volume}{64}
  (\bibinfo{year}{2016}), \bibinfo{pages}{1388--1405}.
\newblock
\showISSN{0030-364X}


\bibitem[(EIA)(2021)]%
        {eia2021monthly}
\bibfield{author}{\bibinfo{person}{U.S. Energy Information~Administration
  (EIA)}.} \bibinfo{year}{2021}\natexlab{}.
\newblock \bibinfo{booktitle}{\emph{June 2021 Monthly Energy Review}}.
\newblock \bibinfo{type}{{T}echnical {R}eport}.
\newblock


\bibitem[Erdo{\u{g}}an and Miller-Hooks(2012)]%
        {erdougan2012green}
\bibfield{author}{\bibinfo{person}{Sevgi Erdo{\u{g}}an} {and}
  \bibinfo{person}{Elise Miller-Hooks}.} \bibinfo{year}{2012}\natexlab{}.
\newblock \showarticletitle{A green vehicle routing problem}.
\newblock \bibinfo{journal}{\emph{Transportation research part E: logistics and
  transportation review}} \bibinfo{volume}{48}, \bibinfo{number}{1}
  (\bibinfo{year}{2012}), \bibinfo{pages}{100--114}.
\newblock


\bibitem[Farr et~al\mbox{.}(2007)]%
        {farr2007shuttle}
\bibfield{author}{\bibinfo{person}{Tom~G Farr}, \bibinfo{person}{Paul~A Rosen},
  \bibinfo{person}{Edward Caro}, \bibinfo{person}{Robert Crippen},
  \bibinfo{person}{Riley Duren}, \bibinfo{person}{Scott Hensley},
  \bibinfo{person}{Michael Kobrick}, \bibinfo{person}{Mimi Paller},
  \bibinfo{person}{Ernesto Rodriguez}, \bibinfo{person}{Ladislav Roth},
  {et~al\mbox{.}}} \bibinfo{year}{2007}\natexlab{}.
\newblock \showarticletitle{The shuttle radar topography mission}.
\newblock \bibinfo{journal}{\emph{Reviews of geophysics}} \bibinfo{volume}{45},
  \bibinfo{number}{2} (\bibinfo{year}{2007}).
\newblock


\bibitem[Fern{\'a}ndez et~al\mbox{.}(2022)]%
        {fernandez2022arc}
\bibfield{author}{\bibinfo{person}{Elena Fern{\'a}ndez},
  \bibinfo{person}{Markus Leitner}, \bibinfo{person}{Ivana Ljubi{\'c}}, {and}
  \bibinfo{person}{Mario Ruthmair}.} \bibinfo{year}{2022}\natexlab{}.
\newblock \showarticletitle{Arc routing with electric vehicles: dynamic
  charging and speed-dependent energy consumption}.
\newblock \bibinfo{journal}{\emph{Transportation Science}}
  (\bibinfo{year}{2022}).
\newblock


\bibitem[Fontana(2013)]%
        {fontanaOptimal2013}
\bibfield{author}{\bibinfo{person}{Matthew~William Fontana}.}
  \bibinfo{year}{2013}\natexlab{}.
\newblock \emph{\bibinfo{title}{Optimal {{Routes}} for {{Electric Vehicles
  Facing Uncertainty}}, {{Congestion}}, and {{Energy Constraints}}}}.
\newblock \bibinfo{thesistype}{Ph.\,D. Dissertation}.
  \bibinfo{school}{Massachusetts Institute of Technology}.
\newblock


\bibitem[Freight(2022)]%
        {uber}
\bibfield{author}{\bibinfo{person}{Uber Freight}.}
  \bibinfo{year}{2022}\natexlab{}.
\newblock \bibinfo{howpublished}{https://www.uber.com/us/en/freight/carrier/}.
\newblock
\newblock
\shownote{Accessed 2022-03-14}.


\bibitem[Goeke and Schneider(2015)]%
        {goekeRouting2015}
\bibfield{author}{\bibinfo{person}{Dominik Goeke} {and}
  \bibinfo{person}{Michael Schneider}.} \bibinfo{year}{2015}\natexlab{}.
\newblock \showarticletitle{Routing a Mixed Fleet of Electric and Conventional
  Vehicles}.
\newblock \bibinfo{journal}{\emph{European Journal of Operational Research}}
  \bibinfo{volume}{245} (\bibinfo{year}{2015}), \bibinfo{pages}{81--99}.
\newblock
\showISSN{0377-2217}


\bibitem[Gonder et~al\mbox{.}(2018)]%
        {gonderFuture2018}
\bibfield{author}{\bibinfo{person}{Jeffrey~D Gonder}, \bibinfo{person}{Aaron~D
  Brooker}, \bibinfo{person}{Eric~W Wood}, {and} \bibinfo{person}{Matthew
  Moniot}.} \bibinfo{year}{2018}\natexlab{}.
\newblock \bibinfo{booktitle}{\emph{Future Automotive Systems Technology
  Simulator (FASTSim) Validation Report}}.
\newblock \bibinfo{type}{{T}echnical {R}eport}. \bibinfo{institution}{National
  Renewable Energy Lab.(NREL), Golden, CO (United States)}.
\newblock


\bibitem[Hancock and Wright(2013)]%
        {hancock2013policy}
\bibfield{author}{\bibinfo{person}{Michael~W Hancock} {and}
  \bibinfo{person}{Bud Wright}.} \bibinfo{year}{2013}\natexlab{}.
\newblock \showarticletitle{A policy on geometric design of highways and
  streets}.
\newblock \bibinfo{journal}{\emph{American Association of State Highway and
  Transportation Officials: Washington, DC, USA}} (\bibinfo{year}{2013}).
\newblock


\bibitem[Hartmann and Funke(2014)]%
        {hartmannEnergyEfficient2014}
\bibfield{author}{\bibinfo{person}{Frederik Hartmann} {and}
  \bibinfo{person}{Stefan Funke}.} \bibinfo{year}{2014}\natexlab{}.
\newblock \showarticletitle{Energy-{{Efficient Routing}}: {{Taking Speed}} into
  {{Account}}}. In \bibinfo{booktitle}{\emph{{{KI}} 2014: {{Advances}} in
  {{Artificial Intelligence}}}} \emph{(\bibinfo{series}{Lecture {{Notes}} in
  {{Computer Science}}})}, \bibfield{editor}{\bibinfo{person}{Carsten Lutz}
  {and} \bibinfo{person}{Michael Thielscher}} (Eds.).
  \bibinfo{publisher}{{Springer International Publishing}},
  \bibinfo{address}{{Cham}}, \bibinfo{pages}{86--97}.
\newblock
\showISBNx{978-3-319-11206-0}


\bibitem[Hassin(1992)]%
        {hassin1992approximation}
\bibfield{author}{\bibinfo{person}{Refael Hassin}.}
  \bibinfo{year}{1992}\natexlab{}.
\newblock \showarticletitle{Approximation schemes for the restricted shortest
  path problem}.
\newblock \bibinfo{journal}{\emph{Mathematics of Operations research}}
  \bibinfo{volume}{17}, \bibinfo{number}{1} (\bibinfo{year}{1992}),
  \bibinfo{pages}{36--42}.
\newblock


\bibitem[Hellstr{\"o}m et~al\mbox{.}(2009)]%
        {hellstrom2009look}
\bibfield{author}{\bibinfo{person}{Erik Hellstr{\"o}m}, \bibinfo{person}{Maria
  Ivarsson}, \bibinfo{person}{Jan {\AA}slund}, {and} \bibinfo{person}{Lars
  Nielsen}.} \bibinfo{year}{2009}\natexlab{}.
\newblock \showarticletitle{Look-ahead control for heavy trucks to minimize
  trip time and fuel consumption}.
\newblock \bibinfo{journal}{\emph{Control Engineering Practice}}
  \bibinfo{volume}{17}, \bibinfo{number}{2} (\bibinfo{year}{2009}),
  \bibinfo{pages}{245--254}.
\newblock


\bibitem[Hu et~al\mbox{.}(2020)]%
        {hu2020battery}
\bibfield{author}{\bibinfo{person}{Xiaosong Hu}, \bibinfo{person}{Le Xu},
  \bibinfo{person}{Xianke Lin}, {and} \bibinfo{person}{Michael Pecht}.}
  \bibinfo{year}{2020}\natexlab{}.
\newblock \showarticletitle{Battery lifetime prognostics}.
\newblock \bibinfo{journal}{\emph{Joule}} \bibinfo{volume}{4},
  \bibinfo{number}{2} (\bibinfo{year}{2020}), \bibinfo{pages}{310--346}.
\newblock


\bibitem[Hwang et~al\mbox{.}(2016)]%
        {hwang2016freight}
\bibfield{author}{\bibinfo{person}{Ho-Ling Hwang}, \bibinfo{person}{Stephanie
  Hargrove}, \bibinfo{person}{Shih-Miao Chin}, \bibinfo{person}{Daniel~W
  Wilson}, \bibinfo{person}{Hyeonsup Lim}, \bibinfo{person}{Jiaoli Chen},
  \bibinfo{person}{Rob Taylor}, \bibinfo{person}{Bruce Peterson}, {and}
  \bibinfo{person}{Diane Davidson}.} \bibinfo{year}{2016}\natexlab{}.
\newblock \bibinfo{booktitle}{\emph{The Freight Analysis Framework Verson 4
  (FAF4)-Building the FAF4 Regional Database: Data Sources and Estimation
  Methodologies}}.
\newblock \bibinfo{type}{{T}echnical {R}eport}. \bibinfo{institution}{Oak Ridge
  National Lab.(ORNL), Oak Ridge, TN (United States)}.
\newblock


\bibitem[Juttner et~al\mbox{.}(2001)]%
        {juttner2001lagrange}
\bibfield{author}{\bibinfo{person}{Alpar Juttner}, \bibinfo{person}{Balazs
  Szviatovski}, \bibinfo{person}{Ildik{\'o} M{\'e}cs}, {and}
  \bibinfo{person}{Zsolt Rajk{\'o}}.} \bibinfo{year}{2001}\natexlab{}.
\newblock \showarticletitle{Lagrange relaxation based method for the QoS
  routing problem}. In \bibinfo{booktitle}{\emph{Proceedings IEEE INFOCOM 2001.
  Conference on Computer Communications. Twentieth Annual Joint Conference of
  the IEEE Computer and Communications Society (Cat. No. 01CH37213)}},
  Vol.~\bibinfo{volume}{2}. IEEE, \bibinfo{pages}{859--868}.
\newblock


\bibitem[Kenworth(2021)]%
        {kenworthT800}
\bibfield{author}{\bibinfo{person}{Kenworth}.} \bibinfo{year}{2021}\natexlab{}.
\newblock \bibinfo{title}{Kenworth T800}.
\newblock
\newblock
\urldef\tempurl%
\url{https://www.kenworth.com/trucks/t800}
\showURL{%
\tempurl}
\newblock
\shownote{Accessed 2021-05-25}.


\bibitem[Liu et~al\mbox{.}(2018)]%
        {liuEnergyEfficient2018}
\bibfield{author}{\bibinfo{person}{Qingyu Liu}, \bibinfo{person}{Haibo Zeng},
  {and} \bibinfo{person}{Minghua Chen}.} \bibinfo{year}{2018}\natexlab{}.
\newblock \showarticletitle{Energy-{{Efficient Timely Truck Transportation}}
  for {{Geographically-Dispersed Tasks}}}. In
  \bibinfo{booktitle}{\emph{Proceedings of the {{Ninth International
  Conference}} on {{Future Energy Systems}}}}. \bibinfo{publisher}{{ACM}},
  \bibinfo{address}{{Karlsruhe Germany}}, \bibinfo{pages}{324--339}.
\newblock
\showISBNx{978-1-4503-5767-8}


\bibitem[Liu et~al\mbox{.}(2017)]%
        {liuConstrained2017}
\bibfield{author}{\bibinfo{person}{Yaqiong Liu}, \bibinfo{person}{Hock~Soon
  Seah}, {and} \bibinfo{person}{Guochu Shou}.} \bibinfo{year}{2017}\natexlab{}.
\newblock \showarticletitle{Constrained Energy-Efficient Routing in Time-Aware
  Road Networks}.
\newblock \bibinfo{journal}{\emph{GeoInformatica}}  \bibinfo{volume}{21}
  (\bibinfo{year}{2017}), \bibinfo{pages}{89--117}.
\newblock
\showISSN{1573-7624}


\bibitem[Lorenz and Raz(2001)]%
        {lorenz2001simple}
\bibfield{author}{\bibinfo{person}{Dean~H Lorenz} {and} \bibinfo{person}{Danny
  Raz}.} \bibinfo{year}{2001}\natexlab{}.
\newblock \showarticletitle{A simple efficient approximation scheme for the
  restricted shortest path problem}.
\newblock \bibinfo{journal}{\emph{Operations Research Letters}}
  \bibinfo{volume}{28}, \bibinfo{number}{5} (\bibinfo{year}{2001}),
  \bibinfo{pages}{213--219}.
\newblock


\bibitem[Maji et~al\mbox{.}(2022)]%
        {majiDACF2022}
\bibfield{author}{\bibinfo{person}{Diptyaroop Maji}, \bibinfo{person}{Ramesh~K.
  Sitaraman}, {and} \bibinfo{person}{Prashant Shenoy}.}
  \bibinfo{year}{2022}\natexlab{}.
\newblock \showarticletitle{{DACF}: day-ahead carbon intensity forecasting of
  power grids using machine learning}. In \bibinfo{booktitle}{\emph{Proceedings
  of the Thirteenth {ACM} International Conference on Future Energy Systems}}
  (Virtual Event, 2022-06-28). \bibinfo{publisher}{{ACM}},
  \bibinfo{pages}{188--192}.
\newblock
\showISBNx{978-1-4503-9397-3}
\urldef\tempurl%
\url{https://doi.org/10.1145/3538637.3538849}
\showDOI{\tempurl}


\bibitem[Maps(2022)]%
        {heremap}
\bibfield{author}{\bibinfo{person}{HERE Maps}.}
  \bibinfo{year}{2022}\natexlab{}.
\newblock \bibinfo{booktitle}{\emph{{Traffic flow using corridor in HERE maps}.
  \url{https://developer.here.com/api-explorer/rest/traffic/flow-using-corridor}}}.
\newblock


\bibitem[Montoya et~al\mbox{.}(2017)]%
        {montoyaElectric2017}
\bibfield{author}{\bibinfo{person}{Alejandro Montoya},
  \bibinfo{person}{Christelle Gu{\'e}ret}, \bibinfo{person}{Jorge~E. Mendoza},
  {and} \bibinfo{person}{Juan~G. Villegas}.} \bibinfo{year}{2017}\natexlab{}.
\newblock \showarticletitle{The Electric Vehicle Routing Problem with Nonlinear
  Charging Function}.
\newblock \bibinfo{journal}{\emph{Transportation Research Part B:
  Methodological}}  \bibinfo{volume}{103} (\bibinfo{year}{2017}),
  \bibinfo{pages}{87--110}.
\newblock
\showISSN{0191-2615}


\bibitem[{OpenStreetMap contributors}(2017)]%
        {OpenStreetMap}
\bibfield{author}{\bibinfo{person}{{OpenStreetMap contributors}}.}
  \bibinfo{year}{2017}\natexlab{}.
\newblock \bibinfo{title}{{Planet dump retrieved from https://planet.osm.org
  }}.
\newblock \bibinfo{howpublished}{\url{ https://www.openstreetmap.org }}.
\newblock


\bibitem[Pelletier et~al\mbox{.}(2017)]%
        {pelletier2017battery}
\bibfield{author}{\bibinfo{person}{Samuel Pelletier}, \bibinfo{person}{Ola
  Jabali}, \bibinfo{person}{Gilbert Laporte}, {and} \bibinfo{person}{Marco
  Veneroni}.} \bibinfo{year}{2017}\natexlab{}.
\newblock \showarticletitle{Battery degradation and behaviour for electric
  vehicles: Review and numerical analyses of several models}.
\newblock \bibinfo{journal}{\emph{Transportation Research Part B:
  Methodological}}  \bibinfo{volume}{103} (\bibinfo{year}{2017}),
  \bibinfo{pages}{158--187}.
\newblock


\bibitem[Schneider et~al\mbox{.}(2014)]%
        {schneiderElectric2014}
\bibfield{author}{\bibinfo{person}{Michael Schneider}, \bibinfo{person}{Andreas
  Stenger}, {and} \bibinfo{person}{Dominik Goeke}.}
  \bibinfo{year}{2014}\natexlab{}.
\newblock \showarticletitle{The {{Electric Vehicle-Routing Problem}} with
  {{Time Windows}} and {{Recharging Stations}}}.
\newblock \bibinfo{journal}{\emph{Transportation Science}}
  \bibinfo{volume}{48} (\bibinfo{year}{2014}), \bibinfo{pages}{500--520}.
\newblock
\showISSN{0041-1655, 1526-5447}


\bibitem[Storandt(2012)]%
        {storandtQuick2012}
\bibfield{author}{\bibinfo{person}{Sabine Storandt}.}
  \bibinfo{year}{2012}\natexlab{}.
\newblock \showarticletitle{Quick and Energy-Efficient Routes: Computing
  Constrained Shortest Paths for Electric Vehicles}. In
  \bibinfo{booktitle}{\emph{Proceedings of the 5th {{ACM SIGSPATIAL
  International Workshop}} on {{Computational Transportation Science}}}}
  \emph{(\bibinfo{series}{{{IWCTS}} '12})}. \bibinfo{publisher}{{Association
  for Computing Machinery}}, \bibinfo{address}{{New York, NY, USA}},
  \bibinfo{pages}{20--25}.
\newblock
\showISBNx{978-1-4503-1693-4}


\bibitem[Strehler et~al\mbox{.}(2017)]%
        {strehler2017energy}
\bibfield{author}{\bibinfo{person}{Martin Strehler}, \bibinfo{person}{S{\"o}ren
  Merting}, {and} \bibinfo{person}{Christian Schwan}.}
  \bibinfo{year}{2017}\natexlab{}.
\newblock \showarticletitle{Energy-efficient shortest routes for electric and
  hybrid vehicles}.
\newblock \bibinfo{journal}{\emph{Transportation Research Part B:
  Methodological}}  \bibinfo{volume}{103} (\bibinfo{year}{2017}),
  \bibinfo{pages}{111--135}.
\newblock


\bibitem[Su et~al\mbox{.}(2021)]%
        {su2021Energy}
\bibfield{author}{\bibinfo{person}{Junyan Su}, \bibinfo{person}{Minghua Chen},
  {and} \bibinfo{person}{Haibo Zeng}.} \bibinfo{year}{2021}\natexlab{}.
\newblock \showarticletitle{Energy Efficient Timely Transportation: A
  Comparative Study of Internal Combustion Trucks and Electric Trucks}. In
  \bibinfo{booktitle}{\emph{Proceedings of the 8th ACM International Conference
  on Systems for Energy-Efficient Buildings, Cities, and Transportation}}
  (Coimbra, Portugal) \emph{(\bibinfo{series}{BuildSys '21})}.
  \bibinfo{publisher}{Association for Computing Machinery},
  \bibinfo{address}{New York, NY, USA}, \bibinfo{pages}{224–225}.
\newblock
\showISBNx{9781450391146}
\urldef\tempurl%
\url{https://doi.org/10.1145/3486611.3492228}
\showDOI{\tempurl}


\bibitem[Teresco(2022)]%
        {metalproject}
\bibfield{author}{\bibinfo{person}{James Teresco}.}
  \bibinfo{year}{2022}\natexlab{}.
\newblock \bibinfo{booktitle}{\emph{{Map-based Educational Tools for Algorithm
  Learning (METAL) project}. \url{https://courses.teresco.org/metal/}}}.
\newblock


\bibitem[Tesla(2022)]%
        {teslaSemi}
\bibfield{author}{\bibinfo{person}{Tesla}.} \bibinfo{year}{2022}\natexlab{}.
\newblock \bibinfo{booktitle}{\emph{Tesla Semi Official Website}}.
\newblock
\newblock
\shownote{Accessed 2022-03-14}.


\bibitem[Wang et~al\mbox{.}(2013)]%
        {wangContextAware2013}
\bibfield{author}{\bibinfo{person}{Yan Wang}, \bibinfo{person}{Jianmin Jiang},
  {and} \bibinfo{person}{Tingting Mu}.} \bibinfo{year}{2013}\natexlab{}.
\newblock \showarticletitle{Context-{{Aware}} and {{Energy-Driven Route
  Optimization}} for {{Fully Electric Vehicles}} via {{Crowdsourcing}}}.
\newblock \bibinfo{journal}{\emph{IEEE Transactions on Intelligent
  Transportation Systems}}  \bibinfo{volume}{14} (\bibinfo{year}{2013}),
  \bibinfo{pages}{1331--1345}.
\newblock
\showISSN{1558-0016}


\bibitem[Xu et~al\mbox{.}(2019)]%
        {xuRide2019}
\bibfield{author}{\bibinfo{person}{Wenjie Xu}, \bibinfo{person}{Qingyu Liu},
  \bibinfo{person}{Minghua Chen}, {and} \bibinfo{person}{Haibo Zeng}.}
  \bibinfo{year}{2019}\natexlab{}.
\newblock \showarticletitle{Ride the {{Tide}} of {{Traffic Conditions}}:
  {{Opportunistic Driving Improves Energy Efficiency}} of {{Timely Truck
  Transportation}}}. In \bibinfo{booktitle}{\emph{Proceedings of the 6th {{ACM
  International Conference}} on {{Systems}} for {{Energy-Efficient Buildings}},
  {{Cities}}, and {{Transportation}}}} \emph{(\bibinfo{series}{{{BuildSys}}
  '19})}. \bibinfo{publisher}{{Association for Computing Machinery}},
  \bibinfo{address}{{New York, NY, USA}}, \bibinfo{pages}{169--178}.
\newblock
\showISBNx{978-1-4503-7005-9}


\bibitem[Zhang et~al\mbox{.}(2022)]%
        {zhang2022optimal}
\bibfield{author}{\bibinfo{person}{Yongzhi Zhang}, \bibinfo{person}{Xiaobo Qu},
  {and} \bibinfo{person}{Lang Tong}.} \bibinfo{year}{2022}\natexlab{}.
\newblock \showarticletitle{Optimal Eco-driving Control of Autonomous and
  Electric Trucks in Adaptation to Highway Topography: Energy Minimization and
  Battery Life Extension}.
\newblock \bibinfo{journal}{\emph{IEEE Transactions on Transportation
  Electrification}} \bibinfo{volume}{8}, \bibinfo{number}{2}
  (\bibinfo{year}{2022}), \bibinfo{pages}{2149--2163}.
\newblock


\bibitem[Zhou et~al\mbox{.}(2020)]%
        {zhou2020minimizing}
\bibfield{author}{\bibinfo{person}{Runzhi Zhou}, \bibinfo{person}{Qingyu Liu},
  \bibinfo{person}{Wenjie Xu}, \bibinfo{person}{Minghua Chen}, {and}
  \bibinfo{person}{Haibo Zeng}.} \bibinfo{year}{2020}\natexlab{}.
\newblock \showarticletitle{Minimizing Emission for Timely Truck Transportation
  with Adaptive Fuel Injection}. In \bibinfo{booktitle}{\emph{Proceedings of
  the 7th ACM International Conference on Systems for Energy-Efficient
  Buildings, Cities, and Transportation}}. \bibinfo{pages}{240--249}.
\newblock


\end{thebibliography}
\bibliographystyle{ACM-Reference-Format}

\appendix
\section{Appendix}

\subsection{An Illustration of the Carbon Footprint Objective}

\begin{figure}[H]
    \includegraphics[width=0.80\linewidth]{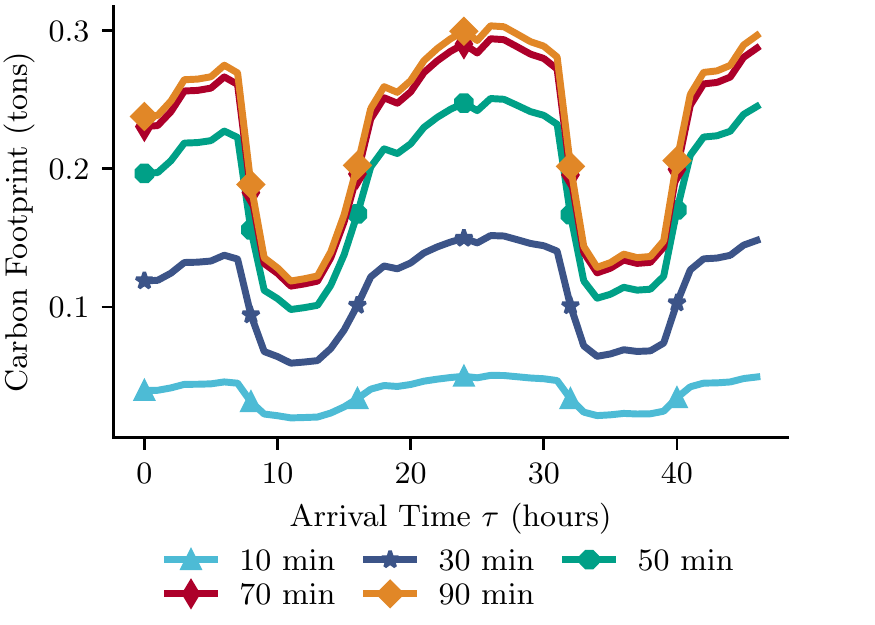}    
    \vspace{-4mm}
    \caption{An illustration of the carbon footprint $F_{v}$ with three variables in California on Jan. 1st, 2022. 
    The carbon footprint is with respect to the arrival time starting at 00:00 AM, Jan. 1st, 2022. 
    Different lines are for different charging times.
    }
    \label{fig:obj}
\end{figure}

\subsection{Discussion on Alternative Approaches}\label{sec:novelty}

In the first alternative,
we directly formulate the CFO problem as a Mixed Integer Programming (MIP). Then we apply an general-purpose method for MIP (e.g., branch and bound) to the formulated problem. As a general observation, the first alternative does not explore the problem structure and only works for problems with hundreds of variables. It usually incurs significant computational burden in larger problems.

The second alternative involves constructing a time-expanded and battery-expanded graph similar in~\cite{strehler2017energy}. For each node $v\in \mathcal{V}$, we add some copies $v_{\beta,\tau}$ of $v$ to the expanded graph. Each copy $v_{\beta,\tau}$ represents a specific SoC $\beta$ and arrival time $\tau$. We add an new edge between $u_{\beta_1, \tau_1}$ and $v_{\beta_2, \tau_2}$ if there is an edge between $(u,v)$ on the original graph and there is a feasible speed choice that corresponds to energy consumption $\beta_2 - \beta_1$ and travel time $\tau_2 - \tau_1$. 
If $v\in \mathcal{V}_c$ is a charging station, we then add an edge between $v_{\beta_1,\tau_1}$ and $v_{\beta_2,\tau_2}$ that corresponds to certain charging decisions. We also assign corresponding weight to this type of edge accordingly. Then the CFO problem can be transformed to a shortest path problem and solved via standard shortest path algorithms e.g., Dijkstra algorithm. 
Although the formulated problem is simple, the second alternative has the drawback that it may exclude feasible solutions of the original problem for inaccurate discretization and may even disconnect origin and destination~\cite{strehler2017energy}. 
Therefore, the constructed expanded graph can be significantly larger than the original graph, and becomes intractable for long-haul tasks. 


The third alternative combines the two directions together. 
We use the battery-expanded graph to deal with the nonlinearity of the energy consumption function and the charging function. If the carbon-intensity function is piecewise linear, then the CFO problem can be formulated as an Mixed Integer Linear Program (MILP) which is less complex than the MIP in the first alternative. We can apply a branch and cut method to the MILP and obtain solutions.
Similar idea can be found in~\cite{fernandez2022arc}.
Similar to the first alternative, the third alternative resorts to solve a MILP with a general-purpose method. It still incurs large computational burden as the problem size grows because MILP is challenging to solve.

Overall, the conceivable alternatives have high model complexity or fail to explore the problem structure, leading to overall high time complexity. In contrast, our novel problem formulation in~\eqref{prob:cfo:a} keeps a low-complexity model and reveals an elegant problem structure that allows us design an dual-subgradient algorithm that has good performance theoretically and empirically.

\begin{table}[!tb]
    \centering
    \begin{tabularx}{\linewidth}{ll}
        \toprule
        Notation & Definition \\
        \midrule
        \midrule
        $\mathcal{G} $ & The original transportation graph \\
        $\mathcal{V}_c $ & The set of charging stations \\
        $\tilde{\mathcal{V}}_c $ & The set of charging stations and source and destination  \\
        $\mathcal{G}_s$ & The stage-expanded graph \\
        $\mathcal{G}_{ex}$ & The extended graph for solving the ILP~\eqref{prob:ilp} \\
        \midrule 
        \midrule 
        $s,d$ & The source and destination respectively \\
        $\beta_0$ & The initial SoC \\
        $T$  & The deadline constraint \\
        $B$ & The battery capacity \\
        $N$ & The maximum number of charging stops\\
        \midrule
        \midrule
        $c_e(t)$ & The energy consumption at road segment $e$  \\
        $\phi_v(t_c, \beta_v)$ & The SoC increment function \\
        $\pi_v(\tau)$ & The carbon intensity at node $v$ \\
        $F_v(\beta, t_c, \tau)$ & The carbon footprint objective \\
        \midrule
        \midrule
        $t_e^i$ & \multicellleft{The travel time at road segment $e$ at stage $i$} \\
        $t_w^{i,v}$ & \multicellleft{The waiting time at charging station $v$ at stage $i$ } \\
        $\beta_v^i$ & The arrival SoC at charing station $v$ at stage $i$  \\
        $\tau_v^i$ & The arrival time at charing station $v$ at stage $i$  \\
        \midrule
        $x_e^i$ & \multicellleft{The binary variable that selects  a road segment $e$ \\ at stage $i$} \\
        \midrule
        $y_v^i$ & \multicellleft{The binary variable that selects a charging station $v$ \\ at stage $i$} \\
        \bottomrule
    \end{tabularx}
    \caption{Key notations.}\label{tab:notation}
\end{table}

\subsection{Proof of Lemma~\ref{lem:alpha:feasible}}\label{sec:lem:alpha:feasible}

\begin{proof}
We denote the initial SoC as $\beta_0$. The SoC when entering the charging stop is no less than $\alpha B$ means that 
\begin{equation}
    \beta_0 - \sum_{i=1}^{n} c_i \geq \alpha B.
\end{equation}
For any $1\leq j < k \leq n$, we have
\begin{align*}
    \sum_{i=j}^{k} c_i &\leq \sum_{i=1}^{n} |c_i| 
    \stackrel{\eqref{ineq:recharge}}{\leq} \frac{1}{1-\alpha} \sum_{i=1}^{n} c_i \\
    &\leq \frac{1}{1-\alpha} \left( \beta_0 - \alpha B \right)  \\
    & \leq \frac{1}{1-\alpha} \left( \beta_0 - \alpha \beta_0 \right) = \beta_0
\end{align*}
Therefore, we have
\begin{align}
    \beta_0 - \sum_{i=j}^{k} c_i \geq 0
\end{align}
for all $1\leq j < k \leq n$. And it directly implies the non-negative SoC on each road segment.
\end{proof}

\begin{figure}[H]
    \includegraphics*[width=.7\linewidth]{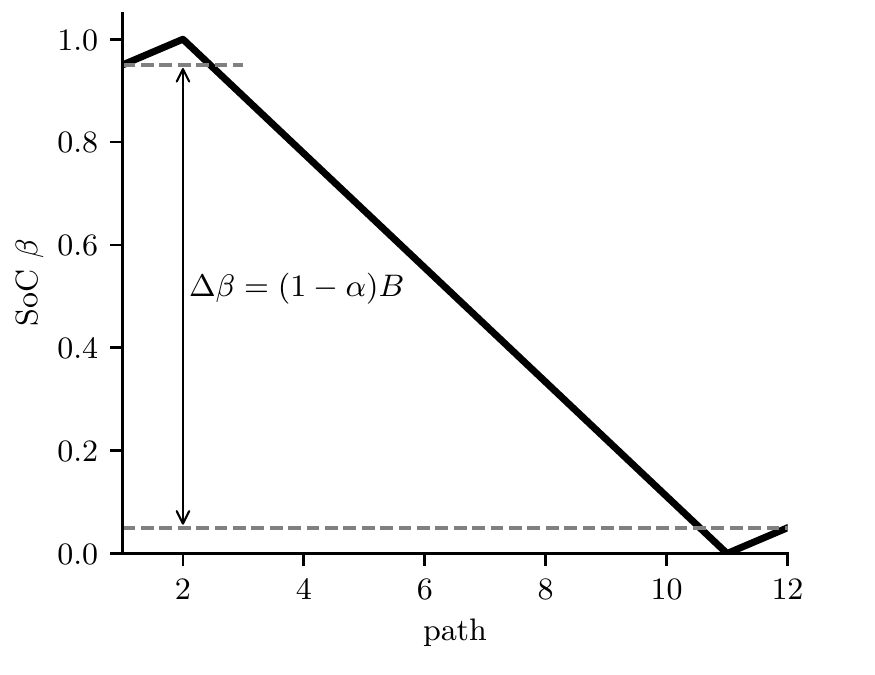}
    \caption{An geometrical interpretation for Lemma~\ref{lem:alpha:feasible}. }
\end{figure}

\subsection{Proof of Theorem~\ref{thm:nphard}}\label{sec:thm:nphard}

\begin{proof}
    We prove \textsf{CFO} is NP-hard by showing that it covers the NP-hard problem \textsf{PASO}~\cite{dengEnergyefficient2017} as a special case. Consider the \textsf{CFO} instance with following parameters: i) the electric truck has a battery with infinite size, i.e. $B=\infty$, ii) there is only one charging station at the origin $s$ with negligible charging time, iii) the battery is empty at the origin, i.e. $\beta_s = \beta_0 = 0$, iv) the carbon-intensity function is constant and equals one at the charging station, i.e. $\pi(\tau) \equiv 1$.
    To solve this instance, one has to charge at the origin with minimum amount of electricity. Then this instance becomes the \textsf{PASO} problem which is known to be NP-hard.
\end{proof}

\subsection{Proof of Theorem~\ref{thm:converge}}\label{sec:thm:converge}

\begin{proof}
    We follow the proof from~\cite{boyd2003subgradient}. We denote the optimal dual variable in~\eqref{prob:dual0} by $\vec{\lambda}^*$ and its norm by $R=\|\vec{\lambda}^*\|$. At $k$-th iteration, we denote by $D_k$ the dual objective and $\vec{\delta}[k]$ its subgradient.  We note that the subgradient is bounded: $\|\vec{\delta}[k]\| \leq H$.
    From the dual update rule~\eqref{eq:dualupdate}, we have
    \begin{subequations}
        \begin{align}
            & \| \vec{\lambda}[k+1]- \vec{\lambda}^* \|^2 
            = \| \vec{\lambda}[k] + \theta_k \vec{\delta}[k] - \vec{\lambda}^* \|^2 \\
            =& \| \vec{\lambda}[k] - \vec{\lambda}^* \|^2 + 2 \theta_k \vec{\delta}[k] \left( \vec{\lambda}[k] - \vec{\lambda}^* \right) + \theta_k^2 \| \vec{\delta}[k] \|^2 \\
            \leq & \| \vec{\lambda}[k] - \vec{\lambda}^* \|^2 + 
            2 \theta_k \left( D_k - D^* \right) + \theta_k^2 H^2 
        \end{align}
    \end{subequations}
    where the inequality follows from the property of the subgradient and the boundedness of the subgradients. By recursively applying the above inequality, we have:
    \begin{subequations}
    \begin{align}
        0 &\leq \| \vec{\lambda}[K+1]- \vec{\lambda}^* \|^2   \\
        & \leq \| \vec{\lambda}[1] -\lambda^* \|^2 + 2 \sum_{k=1}^{K} \theta_k \left( D_k - D^* \right) + K \theta_k^2 H^2 \\
        &= R^2+ \frac{2}{\sqrt{K}} \sum_{k=1}^{K}\left( D_k - D^* \right)  + H^2 \label{eq:thm:converge:1} \\
        &\leq R^2+ \frac{2K}{\sqrt{K}} \left( \overline{D}_K - D^* \right)  + H^2 \label{ineq:them:converge:2}
    \end{align}
    \end{subequations}
    where the equality~\eqref{eq:thm:converge:1} comes from the fact that $\vec{\lambda}[1]=0$, $\|\vec{\lambda}^*\|=R$ and the constant step size $\theta_k = \frac{1}{\sqrt{K}}$. The inequality~\eqref{ineq:them:converge:2} comes from the definition of $\overline{D}_K$. By rearranging the above inequality, we have
    \begin{align}
        D^* - \overline{D}_K \leq C \frac{1}{\sqrt{K}}
    \end{align}
    with constant $C = \frac{R^2+H^2}{2}$.
\end{proof}

\subsection{Proof of Theorem~\ref{thm:optimal}}\label{sec:thm:optimal}

\begin{proof}
    Suppose Alg.~\ref{alg:dualsub} produces a feasible solution $\vec{z}_k$ updated at iteration $k$. We denotes its objective by $ALG$. By weak duality, we have 
    \begin{align*}
        OPT \geq D(\vec{\lambda}[k]) 
        = ALG + \sum_{i=1}^{N+1} \left( \lambda_i^\beta[k] \delta_i^\beta[k] + \lambda_i^\tau[k] \delta_i^\tau[k] \right).
    \end{align*}
    We thus complete the proof by rearranging terms. 
\end{proof}

\subsection{Proof of Proposition~\ref{prop:complexity}}\label{sec:prop:complexity}

\begin{proof}
    The main computational burden of Algorithm~\ref{alg:dualsub} comes from solving the problem~\eqref{prob:dual} in line~\ref{alg:dualsub:primal}. Its complexity has two parts: 

    \trianglebullet~Solve the subproblems~\eqref{prob:w} with complexity $O\left( \left( N+1 \right) |\mathcal{E}| \right)$ and solve the subproblems~\eqref{prob:sigma} with complexity $O\left( N |\mathcal{V}_c| \frac{M^4}{\epsilon_1^4} \right)$. Recall that $M=\max\left\{ t_c^{ub}, t_w^{ub}, B, T\right\}$ is the diameter of the box constraint in subproblem~\eqref{prob:sigma} and $\epsilon_1$ is the accuracy for solving the subproblems~\eqref{prob:sigma}.

    \trianglebullet~For solving the outer problem~\eqref{prob:ilp}, we first need to compute the shortest path for each charging station pair. We apply the single-source all-destination Bellman-Ford~\cite{cormen2022introduction} algorithm for each source in $\mathcal{V}_c$ and for each weight $w_e^i, \forall i \in \left\{ 1,...,N+1 \right\}$. 
    Therefore, the time complexity of this step is 
    {$$O\Bigl( \left( N+1 \right) |\mathcal{V}_c| |\mathcal{V}| |\mathcal{E}|  \Bigr) + O\left( N |\mathcal{V}_c|^2 \right)$$}
    where the second term is the complexity { for the shortest path algorithm on the extended graph of the charge stations}.

    In summary, after rearranging terms, the overall time complexity per iteration of Alg.~\ref{alg:dualsub} is given by
    $$O  \left( \left( N+1 \right) |\mathcal{V}_c| |\mathcal{V}| |\mathcal{E}| + N |\mathcal{V}_c| \frac{M^4}{\epsilon_1^4} \right) $$
\end{proof}

\end{document}